\documentclass[11pt,dvips]{article}
\usepackage{latexsym}
\usepackage{amsmath,amsthm,amsfonts,amssymb} 
\usepackage[T1]{fontenc} 
\usepackage[latin1]{inputenc}
\usepackage{aeguill}
\usepackage{url}

\usepackage{enumerate}
\usepackage{mdwlist} 
\usepackage{graphicx}
\usepackage[all]{xy}

\newtheorem{thm}{Theorem}[section]
\newtheorem{thmbis}{Theorem}
\newtheorem*{thm*}{Theorem}
\newtheorem{dfn}[thm]{Definition} 
\newtheorem*{dfn*}{Definition}

\newtheorem{cor}[thm]{Corollary}
\newtheorem*{cor*}{Corollary}

\newtheorem{prop}[thm]{Proposition} 
\newtheorem*{prop*}{Proposition} 
\newtheorem*{properties*}{Properties} 
 
\newtheorem{lem}[thm]{Lemma} 
\newtheorem*{lem*}{Lemma}

\newtheorem*{claim*}{Claim} 
 
\newtheorem*{fact*}{Fact}

\newtheorem*{qst*}{Question}

\newtheorem*{pb*}{Problem}

\theoremstyle{remark}
\newtheorem*{rem*}{Remark}
\newtheorem{rem}[thm]{Remark}
\newtheorem*{example*}{Examples}
\newtheorem{example}[thm]{Example}

\makeatletter
\edef\@tempa#1#2{\def#1{\mathaccent\string"\noexpand\accentclass@#2 }}
\@tempa\rond{017}
\makeatother

\renewcommand{\phi}{\varphi} 
\newcommand{\m} {^{-1}} 
\newcommand{\eps} {\varepsilon}

\newcommand {\ra} {\rightarrow}

\newcommand{\ol}[1]{\overline{#1}}

\newcommand{\dunion}{\sqcup}

\newcommand{\Dunion}{\bigsqcup} 

\newcommand{\ie} {i.e.\ }

\newcommand {\calc} {{\mathcal {C}}}   
\newcommand {\cald} {{\mathcal {D}}}   
\newcommand {\cale} {{\mathcal {E}}}

\newcommand {\calk} {{\mathcal {K}}}

\newcommand {\calz} {{\mathcal {Z}}}

\newcommand {\bbR} {{\mathbb {R}}}

\newcommand {\bbZ} {{\mathbb {Z}}}   

\newcommand{\grp}[1]{\langle #1 \rangle}

\newcommand{\Out} {{\mathrm{Out}}}
\newcommand{\Aut} {{\mathrm{Aut}}}

\newcommand{\inc}{\subset}
\newcommand {\Z} {{\mathbb {Z}}}   
\newcommand{\cc}{^*_{c}}

\begin{document}

\title{Trees of cylinders and canonical splittings}
\author{Vincent Guirardel, Gilbert Levitt}
\date{}

\maketitle

\begin{abstract}  
  Let $T$ be a tree with an action of a finitely generated group $G$.
  Given a suitable equivalence relation  on the set of edge stabilizers of $T$  (such as commensurability,
  co-elementarity in a  relatively hyperbolic group, or commutation in a
  commutative transitive group),
  we define a tree of cylinders $T_c$. This tree 
 only depends on the deformation space of $T$;
 in particular, it is invariant under automorphisms of $G$ if $T$ is a
  JSJ splitting. We  thus obtain $\Out(G)$-invariant cyclic or
  abelian JSJ splittings. Furthermore, $T_c$ has very strong
  compatibility properties (two trees are compatible if they have a
  common refinement).
\end{abstract}

{\small
\tableofcontents
}

\section{Introduction}

In group theory, a JSJ splitting of a group $G$ is a splitting of $G$ (as a graph of groups)
in which one can \emph{read} any splitting of $G$ (in a given class), and which is \emph{maximal} for this property
\cite{Sela_structure,RiSe_JSJ,DuSa_JSJ,FuPa_JSJ,Bo_cut,Gui_JSJ,GL3}.

In general, JSJ splittings are not unique, and there is a whole space of JSJ splittings
called the \emph{JSJ deformation space} \cite{For_uniqueness,GL3}. A \emph{deformation space} 
\cite{For_deformation, GL1}
  is the set of all splittings whose elliptic subgroups are prescribed (one usually also
adds constraints on edge groups).
A typical 
example of a deformation space is Culler-Vogtmann's Outer Space \cite{CuVo_moduli}. Splittings in the same deformation space are related by a finite sequence of simple moves \cite{For_deformation,ClFo_whitehead, GL1}. 

JSJ deformation spaces being canonical, they are endowed with a natural action of $\Out(G)$.
As they are contractible \cite{GL1,Clay_contractibility},   this usually gives homological information about $\Out(G)$ (see e.g.\  \cite{CuVo_moduli}).

However, deformation spaces of splittings, including JSJ deformation spaces, may also have some bad behaviour.
For instance the action of $\Out(G)$ may fail to be cocompact;   the deformation space may even fail to be finite dimensional.

Much more satisfying is the case when one has a \emph{canonical splitting} rather than just a canonical deformation space.
Such a splitting is invariant under automorphisms; in other words, it is a fixed point for the action of $\Out(G)$ on the deformation space containing it. A typical example is the virtually cyclic JSJ splitting of a one-ended hyperbolic group  $G$ constructed by Bowditch \cite{Bo_cut} 
from the topology of the boundary of $G$. 
Having such a splitting gives much more precise information on $\Out(G)$ (see \cite{Lev_automorphisms,Sela_structure}).
\\

The goal of this paper is to introduce a general construction producing
 a canonical splitting (called the \emph{tree of cylinders}) from a   deformation space $\cald$.
Rather than splittings (or graphs of groups), we think in terms of  trees equipped with an action of $G$.
 We always assume that $G$ is finitely generated.

The construction starts with a class $\cale$ of allowed edge stabilizers, endowed with an \emph{admissible} equivalence relation (see Definition \ref{dfn_admissible}).  
The main examples are commensurability, co-elementarity, commutation 
(see Examples A, B, C below, and Examples \ref{sec_deuxbouts} to \ref{sec_ex_def}).  All trees  are assumed to have edge stabilizers in $\cale$.

Given  a tree $T\in\cald$, the equivalence relation on edge stabilizers partitions edges of $T$ into \emph{cylinders}. 
An essential feature of an admissible relation is that cylinders are connected  (they are subtrees of $T$).
By definition, the \emph{tree of cylinders} of $T$ is  the  tree $T_c$ dual to the covering of $T$ by its cylinders
(see Definition \ref{dfn_arbre_cyl}).

\begin{thmbis}  \label{thm_un}
The tree of cylinders $T_c$ depends only on the deformation space $\cald$ containing $T$.

Moreover, the assignment $T\mapsto T_c$ is functorial: 
 any equivariant map $T\ra T'$
  induces a natural cellular map $T_c\ra T'_c$ (mapping an edge to an edge or a vertex).
\end{thmbis}

 We often say that $T_c$ is the \emph{tree of cylinders of the deformation space}.
It is $\Out(G)$-invariant
   if $\cald$ is. \\

\begin{figure}[htbp]
  \centering
\includegraphics{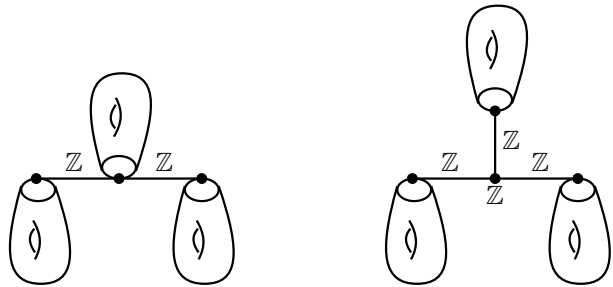}
  \caption{a JSJ splitting of a hyperbolic group and its tree of cylinders}
  \label{fig_ex1}
\end{figure}

\begin{example*}  Consider   the graph of groups pictured on the left of Figure \ref{fig_ex1} 
(edge groups are infinite cyclic and attached to the boundary of punctured tori). Its fundamental group $G$ is hyperbolic, and the splitting depicted is
a cyclic JSJ splitting. Its tree of cylinders is the splitting pictured on the right. It belongs to the same deformation space, but it is $\Out(G)$-invariant; in particular, it has a symmetry of order 3 (it is the splitting constructed in \cite{Bo_cut}).

\begin{figure}[htbp]
  \centering
\includegraphics{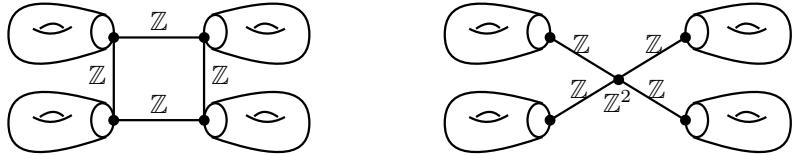}
  \caption{a JSJ splitting of a toral relatively hyperbolic group and its tree of cylinders}
  \label{fig_ex2}
\end{figure}

  Now consider   the graph of groups pictured on the left of Figure \ref{fig_ex2}. Its fundamental group $G$ is
a torsion-free CSA group, actually a toral relative hyperbolic group. The splitting depicted is again
a cyclic JSJ splitting, but   no splitting in its deformation space is $\Out(G)$-invariant.
Its tree of cylinders, which  is $\Out(G)$-invariant, is the splitting pictured on the right. It does not belong to the same deformation space, because of the vertex   with group $\Z^2$.  
\end{example*}

A basic property of $T_c$ is that it is \emph{dominated} by $T$: every   subgroup which is elliptic (\ie fixes a point) in $T$ is elliptic in $T_c$. But, as in the case of Culler-Vogtmann's Outer Space, it may happen that $T_c$ is trivial.
For the construction to be useful, one has to be able to control how far $T_c$ is from $T$, the best situation being when  $T_c$ and $T$  belong to the same deformation space 
(\ie $T_c$ dominates $T$).
One must also control edge stabilizers of  $T_c$, as they   may   fail to be  in $\cale$.

Here are the main examples where this control   is possible.

\begin{itemize}
\item (Example A) $G$ is a hyperbolic group, $\cale$ is the family of two-ended subgroups of $G$, and $\sim$ is the \emph{commensurability} relation.

More generally, 
$G$ is hyperbolic relative to 
subgroups $H_1,\dots H_n$,  the class $\cale$ consists  of the
infinite elementary subgroups,
and $\sim$ is the \emph{co-elementarity} relation; one only considers trees in which each $H_i$ is elliptic   ($A\inc G$ is \emph{elementary} if $A$ is virtually cyclic or contained in a conjugate of some $H_i$, and $A\sim B$ iff $\langle A,B\rangle$ is elementary).
\item (Example B) $G$ is torsion-free and CSA, 
$\cale$ is the class of non-trivial abelian subgroups, and $\sim$ is the \emph{commutation} relation 
(recall that $G$ is CSA if centralizers of non-trivial elements are abelian and malnormal;  $A\sim B$ iff $\langle A,B\rangle$ is abelian).

\item (Example C) $G$ is torsion-free and commutative transitive, 
$\cale$ is the set of infinite cyclic subgroups of $G$, and $\sim$ is the \emph{commensurability} relation. 
\end{itemize}

\begin{thmbis}\label{thm_main}
Let $G$ and $\cale$ be as in Example A (resp.\  B).
Let $T$ be a tree with edge stabilizers in $\cale$, and assume that
non-cyclic abelian subgroups (resp.\ parabolic subgroups)
are elliptic in $T$.
Then: 
\begin{enumerate}
\item $T_c$ has edge stabilizers in $\cale$;
\item $T_c$ belongs to the same deformation space as $T$;
\item $T_c$ is almost $2$-acylindrical.
\end{enumerate}
\end{thmbis}

A tree is \emph{$k$-acylindrical} \cite{Sela_acylindrical} if any segment $I$ of length $> k$ has trivial stabilizer.
 If $G$ has  torsion, we use the notion of \emph{almost acylindricity}:    
the stabilizer of $I$ is finite.

In general, edge stabilizers  of $T_c$ may fail to be in $\cale$. 
In this case, we also
consider the \emph{collapsed tree of cylinders} $T\cc$,
obtained from $T_c$ by collapsing  all  edges whose stabilizer is not in $\cale$.

\begin{thmbis}\label{thm_main2}
Let $G$ and  $\cale$ be as in Example C. 
Let $T$ be a tree with  infinite cyclic edge stabilizers, such that solvable Baumslag-Solitar subgroups of $G$ are elliptic in $T$. 

Then  $T\cc$, $T_c$,  and $T$ lie in the same deformation space, and $T\cc$ is $2$-acylindrical.
\end{thmbis}

Parabolic, abelian, or Baumslag-Solitar subgroups
as they appear in the hypotheses of Theorems \ref{thm_main} and \ref{thm_main2}
 are always elliptic in $T_c$. 
If one does not assume that they are elliptic in $T$, 
they are the only way in which the deformation spaces of $T$ and $T\cc$ differ (see Section \ref{sec_exa} for precise statements).

In general, the deformation space of $T\cc$  may be characterized by the following maximality property: 
 $T\cc$ dominates any tree $T'$ such that $T$ dominates $T'$ and cylinders of $T'$ are bounded  (Proposition  \ref{prop_maxi}). 
In many situations,   
one may replace boundedness of cylinders by acylindricity in the previous  maximality statement (see Section \ref{sec_exa}). 

Our results of Section \ref{sec_exa} may then be interpreted as describing which subgroups must be made elliptic in order to make $T$ acylindrical.  
This is used  in \cite{GL3} to construct (under suitable hypotheses) JSJ splittings of finitely generated groups, using acylindrical accessibility.
\\

Theorems  \ref{thm_main} and \ref{thm_main2} produce a canonical element in the deformation space of $T$. 
In particular, this provides canonical $\Out(G)$-invariant 
JSJ splittings.

\begin{thmbis}\label{thm_can_relh} 
Let $G$ be   hyperbolic relative to $H_1,\dots ,H_n$. 
Assume that $G$ is one-ended relative to $H_1,\dots, H_n$. 
There is an elementary (resp.\ virtually cyclic) JSJ tree relative to $H_1,\dots,H_n$ which is invariant under 
the subgroup of $\Out(G)$ preserving the conjugacy classes of the $H_i$'s. 
\end{thmbis}

 When $G$ is a one-ended hyperbolic group ,   we get   the tree constructed by Bowditch  \cite{Bo_cut}.

\begin{thmbis}\label{thm_can_csa} 
  Let $G$ be a one-ended torsion-free CSA group.  There exists an abelian (resp.\ cyclic) JSJ splitting of $G$ relative to all non-cyclic abelian subgroups, which is  $\Out(G)$-invariant.  
\end{thmbis}

  In particular, one gets canonical cyclic and abelian JSJ splittings of toral relatively hyperbolic groups, including limit groups.
See  \cite{BuKhMi_isomorphism,Pau_theorie,Bo_cut,PaSw_boundaries,ScSw_regular} for other constructions of canonical JSJ splittings.  
As evidenced  by the example  on Figure 2, Theorems \ref{thm_can_relh}  and \ref{thm_can_csa}  do not hold  for   non-relative JSJ splittings.

 Scott and Swarup have constructed in \cite{ScSw_regular} a canonical splitting  over virtually polycyclic groups.  We show in \cite{GL2} that  their splitting coincides (up to subdivision) with 
the tree of cylinders of the   JSJ deformation space   
(with  $\sim$ being commensurability, see Example \ref{sec_commens}).

In a  forthcoming paper \cite{GL4}, we will use   trees of cylinders to study some 
 automorphism groups.
In particular, we will prove that the tree of cylinders of the abelian JSJ  deformation space of
 a CSA group
is the splitting having largest modular automorphism group.
We show that a hyperbolic group (possibly with torsion) has infinite outer automorphism group
if and only if it has a splitting with an infinite order Dehn twist (see also \cite{Carette_automorphism}).
Finally, we prove a structure theorem 
 for automorphisms of one-ended toral relative hyperbolic groups, based on Theorem \ref{thm_can_relh}. 
\\

Another important feature of the tree of cylinders is compatibility.  
A tree $T$ is a \emph{refinement} of $T'$ if $T'$ can be obtained from $T$ by   equivariantly collapsing a set of edges.
Two trees are \emph{compatible} if they have a common refinement.
Refinement is a notion much more rigid than domination: any two trees in Culler-Vogtmann's outer space dominate each other,
but they are compatible if and only if they lie in a common simplex.

\begin{thmbis}\label{thm_comp}
Let $G$ be   finitely presented and one-ended. 
Assume that $\cale$ and $\sim$ are as in Example A, B, or C.

If  $T$ is a JSJ tree  over subgroups of $\cale$,
then $T\cc$ is compatible with every tree with edge stabilizers in $\cale$ (in which parabolic subgroups are elliptic in Example A).
\end{thmbis}

 More general situations are investigated in Section \ref{sec_compat}.   In \cite{GL3}, we use the compatibility properties of $T\cc$ to reprove and generalize the main theorem of \cite{Gui_JSJ}, describing actions of $G$ on $\bbR$-trees.
 
It is also shown in \cite{GL3} that under the hypotheses of Theorem \ref{thm_comp}  the tree  
$T\cc$ is maximal (for domination) among trees which are compatible with every other tree. 
In other words, $T\cc$ belongs to the same deformation space as  the \emph{JSJ compatibility tree} defined in \cite{GL3}.  
\\
 
 The paper is organized as follows.
After preliminaries, we    define admissible equivalence relations and give examples (Section \ref{sec_adm}).
In Section \ref{sec_constr} we define the tree of cylinders. Besides the geometric definition sketched earlier, we give an algebraic one using elliptic subgroups and we show Theorem \ref{thm_un}  (Corollary \ref{cor_deformation}    and Proposition \ref{prop_cy}).  In Section  \ref{sec_espdef} we give basic properties of $T_c$. In particular, we show that most small subgroups of $G$ which are not virtually cyclic are elliptic in $T_c$. We also 
study boundedness of cylinders and acylindricity of $T\cc$,  we show $(T\cc)\cc=T\cc$, and we prove the maximality property of $T\cc$.
Section \ref{sec_exa} gives a description of the tree of cylinders in the main examples and proves Theorems \ref {thm_main} and \ref{thm_main2}.  
 Section \ref{sec_jsj} recalls some material about  JSJ splittings, and proves the existence of canonical JSJ splittings as in 
Theorem \ref{thm_can_relh} and \ref{thm_can_csa}.
In Section \ref{sec_compat} we study compatibility properties of the tree of cylinders and we prove Theorem \ref{thm_comp}.

\section{Preliminaries}\label{sec_prel}
 
In this paper, $G$ will be a fixed finitely generated group. 

Two subgroups $A$ and $B$ are \emph{commensurable} if $A\cap B$ has finite index in both $A$ and $B$. We denote by
$A^g$ the conjugate   $g^{-1}Ag$. The \emph{normalizer} $N(A)$ of $A$ is the set of $g$ such that $A^g=A$. Its
\emph{commensurator}
$Comm(A)$ is the set of $g$ such that $A^g$ is commensurable to $A$. The subgroup $A$ is \emph{malnormal} if $A^g\cap
A\ne\{1\}$ implies
$g\in A$.

 If $\Gamma$ is a graph, we denote by $V(\Gamma)$ its set of vertices and by $E(\Gamma)$ its set of  non-oriented edges. 

A tree  always means a simplicial tree $T$ on which $G$ acts without inversions. Given a family $\cale$ of subgroups of
$G$, an
\emph{$\cale$-tree} is a tree whose edge stabilizers belong to $\cale$.  We denote by $G_v$ or $G_e$ the
stabilizer of a vertex
$v$ or an edge $e$.  

A tree $T$ is \emph{non-trivial} if there is no global fixed point,
\emph{minimal} if there is no proper $G$-invariant subtree.
An element or a subgroup of $G$ is \emph{elliptic} in $T$ if it has a global fixed point. An element   which is not elliptic
is \emph{hyperbolic}. It has an axis   on which it acts as a translation. 

A subgroup $A$ consisting only of elliptic elements fixes a
point if it is finitely generated, a point or an end in general. If a finitely generated  subgroup $A$ is not elliptic, there is 
a unique minimal $A$-invariant subtree.

A group $A$ is \emph{slender} if $A$ and all its subgroups are finitely generated. A slender group acting on a tree fixes a point or leaves a line  invariant (setwise). 

A subgroup $A$ is \emph{small} if it has no non-abelian free subgroups. 
 As in  \cite{GL1}, we could replace smallness by the following weaker property: whenever $G$ acts on  a tree, 
$A$  fixes a point, or an end, or leaves a line  invariant. 

A tree $T$ \emph{dominates} a tree $T'$ if there is an equivariant map $f:T\to T'$. Equivalently, any subgroup which is
elliptic in
$T$ is also elliptic in $T'$. Having the same elliptic subgroups is an equivalence relation on the set of trees, the
equivalence classes are called \emph{deformation spaces} \cite{For_deformation,GL1}.

An equivariant map $f:T\to T'$ between trees \emph{preserves alignment}    if
$x\in[a,b]\implies f(x)\in[f(a),f(b)]$.
 Equivalently, $f$ is a \emph{collapse map}: it is obtained by collapsing certain edges to points. In particular, $f$ does not
fold. 

We   say that
$T$ is a
\emph{refinement} of
$T'$ if there is a collapse map $f:T\to T'$.  Two  trees $T$ and $T'$ are \emph{compatible} if they have a common
refinement.

\section{Admissible relations}\label{sec_adm}

 Let $\cale$ be a class of subgroups of $G$, stable under conjugation. 
It should   not be  stable under taking subgroups (all trees of cylinders are trivial if it is), but it  usually is 
sandwich-closed: if $A\subset H\subset B$ with $A,B\in\cale$, then $H\in \cale$.
 Similarly, $\cale$ is usually invariant under $\Aut(G)$.

\begin{dfn}\label{dfn_admissible}
  An equivalence relation $\sim$ on $\cale$ is \emph{admissible} if the following axioms hold for any $A,B\in\cale$:
  \begin{enumerate}
\item If $A\sim B$, and $g\in G$, then $A^g\sim B^g$.
  \item If   $A\subset B$, then $A\sim B$.
  \item Let $T$ be an $\cale$-tree. If $A\sim B$, and $A$, $B$
fix
$a,b\in T$ respectively, then for each edge
$e\subset [a,b]$ one has
$G_e\sim A\sim B$. 
  \end{enumerate}
\end{dfn}

More generally, one can require 3 to hold only for certain $\cale$-trees.
In particular, given subgroups $H_i$, we say that $\sim$ is \emph{admissible relative to the $H_i$'s} if 3 holds 
for all $\cale$-trees $T$  in
which each $H_i$ is elliptic.

When proving that a relation is admissible, the only non-trivial part usually is   axiom 3.
The following criterion will be useful:
\begin{lem}\label{lem_crit} If  $\sim$ satisfies 1, 2, and   axiom 3' below, then it is admissible.
\begin{enumerate}
\item[3'.] Let $T$ be an $\cale$-tree. If $A\sim B$, and  $A,B $ are elliptic in  
$T $, then $\grp{A,B}$ is also elliptic in $T$.
\end{enumerate}
\end{lem}

\begin{proof} We show axiom 3.
Let $c$ be  a vertex  fixed by $\grp{A,B}$. Since $[a,b]\subset [a,c]\cup [c,b]$, one may assume $e\subset [a,c]$.
One has $A\subset G_e$ because $A$ fixes $[a,c]$,  so that $G_e\sim A$ by axiom 2.
\end{proof}

Given an admissible relation $\sim$ on $\cale$, we shall associate a tree of cylinders  $T_c$ to any $\cale$-tree $T$. If $\sim$ is
only admissible relative to subgroups $H_i$, we require that each $H_i$ be elliptic in $T$.  

Here
are the main examples to which this will apply. 

\newenvironment{subsectionExample}{}{}
\begin{subsectionExample}
\makeatletter \def\@seccntformat#1{Example \csname the#1\endcsname. } \makeatother

\subsection{Two-ended subgroups}\label{sec_deuxbouts}

Let $\cale$ be the set of two-ended subgroups of $G$ (a group $H$ is two-ended if and only if 
 some finite index subgroup of $H$ is infinite cyclic).
Commensurability, defined by $A\sim B$ if $A\cap B$ has finite index in both $A$ and $B$,
is an admissible equivalence relation.

Axiom 1 is clear.  If $A,B$ are two-ended subgroups with $A\subset B$, then $A$ has finite index in $B$, so 2 holds. If
$A,B\in\cale $ are commensurable and fix   points $a,b\in T$, then
$A\cap B$ fixes $[a,b]$. If $e\subset [a,b]$, then $A \cap B\subset G_e$ (with finite index), and $G_e\sim A\cap B\sim A$.

Commensurability is admissible also on the set of infinite cyclic subgroups.

\subsection{Commensurability}\label{sec_commens} 

This class of examples generalizes the previous one.  It is used in \cite{GL2}.

Let $\cale$ be a conjugacy-invariant family of subgroups of $G$ such that:
\begin{itemize}
\item any subgroup $A$ commensurable with some $B\in\cale$ lies in $\cale$;
\item if $A,B\in \cale$ are such that $A\subset B$, then $[B:A]<\infty$.
\end{itemize}

As in the previous example, one easily checks that commensurability  is an admissible equivalence relation on $\cale$. 

For instance, $\cale$ may consist of all subgroups of $G$ which are virtually $\bbZ^n$ for some fixed $n$, or all
subgroups which are virtually polycyclic   of Hirsch length exactly
$n$. The classes of edge groups $\calz\calk$ considered by Dunwoody-Sageev in \cite{DuSa_JSJ} also fit into this example.

\subsection{Co-elementarity 
 (splittings relative to parabolic groups)}\label{sec_hyp_rel}

Let $G$ be hyperbolic relative to   finitely generated subgroups 
$H_1,\dots,H_n$. Recall that a subgroup of $G$ is \emph{parabolic} if it is
contained in a conjugate of an $H_i$,    \emph{elementary} if 
it is finite, or two-ended, or   parabolic. If every $H_i$ is slender (resp.\ small),
 the elementary subgroups are the same as the slender (resp.\ small) 
subgroups.

Let $\cale$ be the class of \emph{infinite} elementary subgroups of $G$.
Let $\sim$ be the co-elementarity relation, defined by saying that $A,B\in \cale$ are \emph{co-elementary} if and only if
$\grp{A,B}$ is elementary.

\begin{lem}
Co-elementarity is an equivalence relation on $\cale$.
\end{lem}

\begin{proof}
Let $X$ be a proper hyperbolic metric space on which $G$ acts properly discontinuously by isometries
(not cocompactly), as in \cite{Bow_relhyp}.  
We claim that two infinite elementary subgroups $A$ and $B$ are co-elementary if and only if they have the same limit
set in
$\partial_\infty(X)$. The transitivity of $\sim$ follows.

For each infinite elementary subgroup $A$, denote by $\Lambda(A)$ its limit set, which consists of one or two points.
If $\Lambda(A)=\Lambda(B)$, then $\grp{A,B}$ stabilizes $\Lambda(A)$ and is elementary.
For the converse, it suffices to prove that $A\subset B$ implies $\Lambda(A)=\Lambda(B)$ whenever $A,B$ are infinite
elementary. Since $A\subset B$, we have $\Lambda(A)\subset\Lambda(B)$.
If $\Lambda(A)=\{\alpha\}\neq \Lambda(B)=\{\alpha,\beta\}$, then $A$
fixes $\alpha$ and $\beta$, so $A$ preserves the geodesics joining $\alpha$ to $\beta$, a contradiction with $\Lambda(A)=\{\alpha\}$.
\end{proof}

\begin{rem} \label{rem_ma}
Note that any infinite elementary subgroup is contained in a unique maximal one, namely the stabilizer of its limit set.
This subgroup is two-ended or conjugate to an $H_i$. 
\end{rem}

\emph{Co-elementarity is admissible relative to the $H_i$'s.}
 Indeed, let us show  axiom 3', assuming that each $H_i$ is elliptic in $T$. The group 
  $\grp{A,B}$ is elementary. It is clearly elliptic if it is parabolic. If it is two-ended, then it contains $A$ with finite index,
so  is elliptic in
$T$ because $A$ is elliptic.  

\subsection{Co-elementarity 
 (arbitrary splittings)}\label{sec_hyp_abs}

We again assume that $G$ is hyperbolic relative to finitely generated subgroups $H_1,\dots, H_n$, and $\cale $ is the family of infinite elementary
subgroups.  If we   want to define a tree of cylinders for any $\cale$-tree $T$, without
assuming that the
$H_i$'s are elliptic in $T$,  we need $\sim$ to be admissible (not just relative to the $H_i$'s).

\begin{lem} \label{lem_arb}
Suppose that  each $H_i$ is   finitely-ended   (\ie $H_i$ is finite, virtually cyclic, or one-ended). Then
co-elementarity is admissible on $\cale$.
\end{lem}

\begin{proof}
 We show axiom 3, with $T$ any $\cale$-tree with infinite elementary edge stabilizers.
Let $H$ be the maximal elementary subgroup containing both $A$ and $B$ (see Remark \ref{rem_ma}).
It is two-ended or  conjugate to an infinite  $H_i$.

If $H$ fixes a point $c\in T$, we argue as in the proof of Lemma \ref{lem_crit}.
If $H$ is two-ended, it contains $A$ with finite index   and therefore   fixes a point in $T$.
Thus, we can assume that $H$ is one-ended and does not fix a point in $T$.  As $H$ is finitely generated, there is  
a minimal
$H$-invariant subtree  
$T_H\subset T$. 

The segment $[a,b]$ is contained in $[a,a']\cup [a',b']\cup [b',b]$ where 
$a'$ (resp.\ $b'$) is the projection of $a$ (resp.\ $b$) on
$T_H$. 
If $e\subset [a,a']$, one has $A\subset G_e$ since $A$ fixes $[a,a']$, so that  $A\sim G_e$.
The same argument applies if $e\subset [b,b']$.
Finally, assume $e\subset [a',b']\subset T_H$. 
Since $H$ is one-ended, the subgroup of $H$ stabilizing $e$ is infinite, so is in $\cale$.
Thus $G_e\sim (G_e\cap H)\sim H\sim A$. Axiom 3 follows.
\end{proof}

\subsection{Commutation 
}\label{sec_ct}

Recall that $G$ is \emph{commutative transitive} if the commutation relation is a transitive relation on 
$G\setminus\{1\}$ (\ie non-trivial elements have abelian centralizers). For example, torsion-free   groups which are
hyperbolic relative to abelian subgroups are commutative transitive.

Let $G$ be a commutative transitive group, and let $\cale$ be the class of its \emph{non-trivial}
 abelian subgroups.

\begin{lem}
  The commutation relation, defined by $A\sim B$ if $\grp{A,B}$ is abelian, is admissible.
\end{lem}
\begin{proof} It is an equivalence relation because of commutative transitivity.
 Axioms 1 and 2 are clear. We show axiom 3' (see Lemma \ref{lem_crit}). 
If   $\grp{A,B}$ does not fix a point in $T$, the fixed point sets of $A$ and $B$ are disjoint, and $\grp{A,B}$ contains a
hyperbolic element $g$. Being abelian, $\grp{A,B}$ then acts by translations on a line (the axis of $g$). But  since $A$ and
$B$ are elliptic,
$\grp{A,B}$   fixes this line pointwise, a contradiction. 
\end{proof}

We will  also consider a more restricted situation.   Note that, in  a commutative transitive
group, any  non-trivial abelian subgroup is contained in  a  unique  maximal abelian subgroup, namely its centralizer.

\begin{dfn}
  $G$ is CSA if it is commutative transitive, and its maximal abelian subgroups are malnormal.
\end{dfn}

Since CSA is a closed property in the space of marked groups,
 $\Gamma$-limit groups for $\Gamma$ torsion-free hyperbolic are CSA (\cite{Sela_diophantine7});
  also see \cite{Champetier_espace} for wilder examples.

\subsection{Finite groups}\label{sec_fini}

Let $q\ge1$ be an integer.
Let $\cale$ be the family of  all  subgroups of $G$ of  cardinality $q$. It is easily checked that equality is an admissible
equivalence relation on
$\cale$. This example will be used in \cite{GL4}.

\subsection{The equivalence relation of a deformation space}\label{sec_ex_def}

\begin{lem}
Let $\cale$ be any conjugacy-invariant family of subgroups of $G$. The equivalence relation generated by inclusion is
admissible relative to $\cale $.
\end{lem}

\begin{proof}
We have to show axiom 3, under the additional hypothesis that all groups of $\cale$ are elliptic in $T$.  
Since $A\sim B$, we can find subgroups  $A_0=A,A_1,\dots ,A_n=B$ in $\cale$
such that  $A_i,A_{i+1}$ are nested (one contains the other).
Let $u_0=a,u_1,\dots,u_n=b$ be  points of $T$ fixed by  $A_0,\dots ,A_n$ respectively.
Let $i$ be such that $e\subset [u_i,u_{i+1}]$, and assume for instance that $A_i\subset A_{i+1}$.
Then $A_i$ fixes $e$, so $A_i\subset G_e$ and $G_e\sim A_i \sim A$.
\end{proof}

In particular, let $\cald$ be a deformation space (or a restricted deformation space in the sense of Definition 3.12
of \cite{GL1}).  The   relation generated by inclusion is
admissible on the family $\cale$ consisting of generalized edge groups of
reduced trees in
$\cald$ (see Section 4 of  \cite{GL1}).

\end{subsectionExample}

\section{The basic construction}\label{sec_constr}

Let $\sim$ be   an admissible equivalence relation on $\cale$. We   now associate a tree of cylinders $T_c$ to any
$\cale$-tree $T$. If $\sim$ is only admissible relative to subgroups $H_i$, we require
that each $H_i$ be elliptic in $T$.
\subsection{Cylinders}

\begin{dfn}
  Let $T$ be an $\cale$-tree  (with
  each $H_i$   elliptic in $T$ in the relative case). 

Define an equivalence relation  $\sim_T$ on the set  of (non-oriented)  
edges of $T$ by: $e\sim _T e'\iff G_e\sim G_{e'}$. 
A \emph{cylinder} of $T$ is an equivalence class $Y$. We identify $Y$ with the union of its edges, 
a subforest  of $T$.
\end{dfn}

A key feature of cylinders is their
connectedness:

\begin{lem}\label{lem_conn}
  Every cylinder is a subtree.
\end{lem}

\begin{proof}
  Assume that $G_e\sim G_{e'}$. By axiom 3, any edge $e''$ contained in the arc joining $e$ to $e'$ satisfies
$G_e\sim G_{e'}\sim G_{e''}$, thus belongs to the same cylinder as $e$ and $e'$. 
\end{proof}

Two distinct cylinders meet in at most one point.
One   can then define the tree of cylinders of $T$ as the tree $T_c$ \emph{dual} to the covering of $T$ by its cylinders, as
in \cite[Definition 4.8]{Gui_limit}:

\begin{dfn}\label{dfn_arbre_cyl}
The   \emph{tree of cylinders} of $T$ is the bipartite tree $T_c$ with vertex set $V(T_c)=V_0(T_c)\dunion V_1(T_c)$
defined as follows:
\begin{enumerate}
\item $V_0(T_c)$ is the set of vertices $x$ of $T$ belonging to (at least) two distinct cylinders;
\item $V_1(T_c)$ is the set of cylinders $Y$ of $T$;
\item there is an edge $\varepsilon =(x,Y)$ between $x\in V_0(T_c)$ and $Y\in V_1(T_c)$ if and only if $x$ (viewed as a
vertex of $T$) belongs to $Y$ (viewed as a subtree of $T$).
\end{enumerate}
\end{dfn}

Alternatively, one can define the  \emph{boundary} $\partial Y$ of a cylinder $Y$ as the set of vertices of $Y$ 
belonging to another cylinder, and obtain
$T_c$   from $T$ by replacing each cylinder by the cone on its
boundary.

It is easy to see that $T_c$ is indeed a  tree \cite{Gui_limit}.
Here are a few other simple observations.  

The group $G$ acts on $T_c$.  It follows from \cite[Lemma 4.9]{Gui_limit} that $T_c$ is minimal if $T$ is minimal. But $T_c$ may be a point, for instance
if all edge stabilizers of $T$ are trivial.

Any vertex  stabilizer
$G_v$ of $T$ fixes a point in $T_c$: the vertex $v$ of $V_0(T_c)$ if $v$ belongs to two cylinders, the vertex $Y$ of
$V_1(T_c)$ if
$Y$ is the only cylinder containing $v$.  In other words, $T$ dominates $T_c$. 

A vertex $x\in V_0(T_c)$ may be viewed either as a vertex of $T$ or as a vertex of $T_c$; its
stabilizer in $T_c$ is the same as in $T$. 
The
stabilizer of a vertex in
$V_1(T_c)$ is the (global) stabilizer
$G_Y$  of a cylinder $Y\subset T$; it may fail to be  elliptic in
$T$ (for instance if $T_c$ is a point and $T$ is not), so $T_c$ does not always dominate $T$. This will be studied in Sections
\ref{sec_espdef} and \ref{sec_exa}. 

   Let us now consider edge stabilizers. We note: 
  \begin{rem} \label{rem_gros}
Edge stabilizers of $T_c$ are elliptic in
$T$, and they    always contain a group in $\cale$: 
if $\eps=(x,Y)$ and   $e$ is an edge of $Y$ incident on $x$, then $G_\eps\supset G_e$.
 \end{rem}
 
However, edge stabilizers of $T_c$ are not necessarily in $\cale$.
For this reason, it is convenient to introduce the \emph{collapsed
tree  of cylinders}: 

\begin{dfn}\label{dfn_T0}
Given an $\cale$-tree $T$, 
the \emph{collapsed tree of cylinders} $T\cc$
is the $\cale$-tree obtained from $T_c$ by collapsing all edges whose stabilizer is not in $\cale$.
\end{dfn}

\subsection{Algebraic interpretation}\label{sec_alg}

We give a more algebraic definition of $T_c$, by viewing it as a subtree of a bipartite graph $Z$
 defined algebraically using only information on $\cale$, $\sim$, and elliptic subgroups of $T$.
This will  make it clear  that $T_c$ only depends on the deformation space of $T$ (Corollary \ref{cor_deformation}).
We motivate the definition of $Z$ by a few observations.

If $Y\subset T$ is a cylinder, all its edges have equivalent stabilizers, and we can associate to $Y$ an equivalence class
$\calc\in\cale/\sim$. We record the following for future reference. 

 \begin{rem}\label{rem_classe}
Given an edge $\varepsilon
=(x,Y)$ of $T_c$, let $e$ be an edge of $Y$ adjacent to $x$ in $T$. Then $G_e$ is a
representative of the class $\calc$ which is contained in $G_\varepsilon
=G_x\cap G_Y$. If
$G_\varepsilon
$ belongs to $\cale$, then it is in $\calc$ by axiom 2 of admissible relations. In
particular, if all edge stabilizers of $T_c$ are  in
$\cale$, then $(T_c)_c=T_c$. Also note that $G_Y$ represents $\calc$ if $G_Y\in\cale$.
\end{rem}

Depending on the context, it may be  convenient to think of a cylinder  either as a set of edges of
$T$, or a subtree $Y$ of
$T$, or a vertex of $T_c$, or   an equivalence class $\calc$.  Similarly,  there are several ways to think of $x\in V_0(T_c)$:
as a vertex of $T$, a vertex of $T_c$, or an elliptic subgroup $G_x$. 

If $x$  is a vertex of $T$ belonging to  two cylinders, then its stabilizer $G_x$
is not   contained in a group of $\cale$: otherwise all edges of $T$ incident to $x$ would have
equivalent stabilizers by axiom 2, and
$x$ would belong to only one cylinder.

More generally, let $v $ be any vertex of $T$ whose stabilizer  is not contained in a group of $\cale$. Then $G_v$
fixes $v$ only,
  and is a maximal elliptic subgroup. Conversely, let $H$ be a subgroup which is elliptic in $T$, is not contained in a
group of
$\cale$, and is maximal for these properties. Then $H$ fixes a unique vertex $v$ and equals $G_v$.

\begin{dfn}\label{dfn_graphe} Given an $\cale$-tree $T$, 
let $Z$ be the bipartite graph   with vertex set $V(Z)=V_0(Z)\dunion V_1(Z)$
defined as follows:
\begin{enumerate}
\item $V_0(Z)$ is the set of subgroups $H$ which are elliptic in $T$, not contained in  a group of $\cale$, and
maximal for these properties;
\item $V_1(Z)$ is the set of equivalence classes $\calc\in\cale/\sim$;
\item there is an edge $\varepsilon $ between $H\in V_0(Z)$ and $\calc\in V_1(Z)$ if and only if  $H$ contains a group of
$\calc$. 
\end{enumerate}
\end{dfn}

As previously observed, $V_0(Z)$ may be viewed as the set of vertices of $T$ whose stabilizer is not contained in a
group of $\cale$.

It also follows from the previous observations that there is a natural embedding  of   bipartite graphs $j:T_c\ra Z$: 
for  $v\in V_0(T_c)$, we define $j(v)=G_v\in V_0(Z)$;
for  $Y\in V_1(T_c)$, with associated  equivalence class $\calc$, we  define $j(Y)=\calc\in V_1(Z)$.
Note that  $j$ is well defined on $E(T_c)$ since  adjacent vertices of $T_c$ have  adjacent 
images in
  $Z$ by Remark \ref{rem_classe}.
  
  A vertex $H\in V_0(Z)$ is a stabilizer $G_v$ for a unique $v\in T$. It is in $j(T_c)$ if and only if  $v$ belongs to two
cylinders. A vertex
$\calc\in V_1(Z)$ is in $j(T_c)$ if and only if some representative of $\calc$ fixes an edge of $T$.

\begin{prop}\label{segmentcinq}
Assume that the action of $G$   on $T$ is minimal and non-trivial, and that $T_c$ is not a point. 
Then $j(T_c)$ is the set of edges and vertices of $Z$ which 
are contained in the central edge of a segment of length 5 of $Z$.
\end{prop}

\begin{proof}
The action of $G$ on $T_c$ is minimal (see above) and non-trivial, so any edge of $T_c$ and of $j(T_c)$ is the central edge of a segment
of length 5.
The converse is  an immediate consequence of items 1, 3, 4 of the following  lemma.
\end{proof}  

\begin{lem}\label{segmentcinqc} 
  \begin{enumerate} 
  \item  $v\in V_1(Z)$ belongs to $ j(T_c)$ if and only if   $v$ has valence at least $2$ in $Z$.
  \label{cl1}
  \item  if $x\in T$ and $G_x\in V_0(Z)$ is adjacent to $j(Y)\in j(T_c)$, then $x\in Y$. \label{cl2}
  \item  an element of $  V_0(Z)$ belongs to $ j(T_c)$   if and only if  it has at least 2 neighbours in $j(T_c)$.\label{cl3}
  \item an edge of $Z$ lies in $j(E(T_c))$ if and only if its endpoints are in $j(V(T_c))$.
  \label{cl4}
  \end{enumerate}
\end{lem}

\begin{proof}
In statements \ref{cl1}, \ref{cl3} and \ref{cl4}, the direct implications are clear.

\ref{cl1}. Consider $\calc \in V_1(Z)$ of valence at least $2$. It has representatives $A\subset G_x$ and $B\subset
G_{x'}$ with $x,x'$ distinct points of $ T$. By axiom 3 of admissible relations, all edges in $[x,x']$ have stabilizer in $\calc$. This shows $\calc\in j(V_1(T_c))$.

\ref{cl2}. Consider $G_x\in V_0(Z)$ adjacent to $\calc=j(Y)\in V_1(Z)$ for some cylinder $Y$ of $T$. Then
$\calc$ has a representative $A\subset G_x$. Let $e$ be an edge of $Y$, so that $A\sim G_e$. 
By axiom 3, the cylinder $Y$ contains the arc joining $e$ to $x$, so $x\in Y$. 

\ref{cl3}. If $G_x\in V_0(Z)$ is adjacent to $j(Y)$ and $j(Y')$, then $x\in Y\cap Y'$ by statement \ref{cl2},
so  $x\in V_0(T_c)$, and $G_x\in j(V_0(T_c))$.

\ref{cl4}. Let $j(x)$ and $j(Y)$ be vertices of $j(T_c)$ joined by an edge of $Z$.
 Then $x\in Y$ by statement \ref{cl2}, so $x$ and $Y$ are joined by an edge in $T_c$.
\end{proof}

\begin{cor}\label{cor_deformation}    
If $T,T' $ are minimal  non-trivial $\cale$-trees belonging to the same deformation space, there is 
a canonical equivariant isomorphism between their trees of cylinders.
\end{cor}

\begin{proof} The key observation is that the graph  $Z$ is defined purely in terms of the elliptic subgroups, which are
the same for
$T$ and $T'$. The corollary is trivial if   $T_c,T'_c$ are both  points. If
$T_c$ is not a point, then $V_1(Z)$ contains   infinitely many vertices of valence $\ge2$.  As shown above (statement 1 of  
 Lemma \ref{segmentcinqc}), this implies
that $T'_c$ is not a point, and the result follows directly from Proposition
\ref{segmentcinq}. 
\end{proof}

\subsection{Functoriality}\label{sec_fonct}
There are at least two other ways of proving   that $T_c$ only depends on the deformation space of $T$.  One is
based on the fact that any two trees in the same deformation space are connected by a finite sequence of elementary
expansions and collapses \cite{For_deformation}. One  checks that these moves do not change $T_c$. 

Another approach is to study the effect of an equivariant map $f:T\to T'$ on the trees of cylinders. We always assume
that $f$ maps a vertex to a vertex, and an edge to a point or an edge path. If $T'$ is minimal, any $f$ is surjective.
We say that $f$ is \emph{cellular} if it maps an edge to an edge or a vertex.

\begin{prop}\label{prop_cy} Let 
 $T,T' $ be minimal $\cale$-trees, and $f:T\ra T'$ an equivariant map.
Let $T_c,T'_c$ be the trees of cylinders of $T$ and $T'$.
Then $f$ induces a   cellular equivariant map $f_c:T_c\ra T'_c$. This map does not depend on  
$f$, and is functorial in the sense that 
 $(f\circ g)_c=f_c\circ g_c$.
\end{prop}

Corollary \ref{cor_deformation} easily follows from this proposition. 
The proposition   may be proved by considering the bipartite graph $Z$, but we give a geometric argument. 

\begin{proof} 
We may assume that $T'_c$ is not a point.

\begin{lem}\label{lem_image_cyl}
Consider an equivariant map $f:T\ra T'$.
For each cylinder $Y\subset T$, the image $f(Y)$ is either a cylinder
$Y'$ of 
$T'$ or a point $p'\in T'$.
\end{lem}

\begin{proof}
If an edge $e'$ of $T'$ is contained in the image of an edge $e\subset Y$, then $G_{e'}$ contains $G_e$,
hence is equivalent to $G_e$ by axiom 2. 
This shows that $f(Y)$, if not a point, is contained in a unique cylinder $Y'$. Conversely, if $e'$ is an edge
of $Y'$, then any edge $e$ such that $e'\subset f(e)$ satisfies $G_e\subset G_{e'}$ hence is in $Y$. Thus $Y'=f(Y)$. 
\end{proof}

We say that a cylinder $Y$ is \emph{collapsed} if $f(Y)$ is a point $p'\in
T'$. We claim that such a $p'$ belongs to two distinct cylinders of $T'$, so represents an element of $ V_0(T'_c)$.
Consider the union of all collapsed cylinders, and    the component
containing
$Y$. It is not the whole of $T$, so by minimality of $T$ it has at least two
boundary points. These points belong to distinct non-collapsed 
cylinders whose images are the required cylinders containing $p'$. 

Also
note that, if $x\in T$ belongs to two cylinders, so does $f(x)$. This is
clear if $x$   belongs to no collapsed cylinder, and follows from the
previous fact if it does.

This allows us to
 define $f_c $ on vertices of $T_c$, by sending  $x\in V_0(T_c)$ to $f(x)\in
V_0(T'_c)$, and $Y\in V_1(T_c)$  to  $Y'\in V_1(T'_c)$ or $p'\in V_0(T'_c)$. If
$(x,Y)$ is an edge of
$T_c$,  then $f_c(x)$ and $f_c(Y)$ are equal or adjacent in $T'_c$.  

We may describe $f_c$ without referring to $f$, as follows. The image of $x\in V_0(T_c)$ 
is the unique point of $T'$ fixed by $G_x$.   The image of
$Y$ is the unique  cylinder whose edge stabilizers are equivalent to those of $Y$,
or the unique point of $T'$ fixed by stabilizers of edges of $Y$.  Functoriality is
easy to check.
\end{proof}

\begin{rem}\label{lem_fc_collapse}  Note that, if two edges $(x_1,Y_1)$ and $(x_2,Y_2)$ of $T_c$ are mapped by $f_c$
onto the same edge of
$T'_c$, then $Y_1=Y_2$. In particular,  if the restriction of $f$ to each cylinder is either constant or injective, then $f_c$
preserves alignment. \end{rem}

\section{General properties  }\label{sec_espdef}

\subsection{The deformation space of $T_c$}

We fix $\cale$ and an admissible relation $\sim$.
We have seen that $T$ always dominates $T_c$: any vertex stabilizer of $T$ fixes a point in $T_c$. Conversely,   $T_c$ has two types of vertex stabilizers. 
If $v\in V_0(T_c)$, then its stabilizer is a vertex stabilizer of $T$,
and $G_v\notin \cale$ (see  Subsection \ref{sec_alg}). On the other hand, the
stabilizer  $G_Y$ of a vertex  
$Y\in V_1(T_c)$   may fail to be  elliptic in
$T$. This means that $T_c$ is not necessarily in the same deformation
space as
$T$.

There are various ways to think of   $G_Y$. It consists of those $g\in G$ that map the cylinder  
$Y$ to itself. If
$e$ is any edge in $Y$, then $G_Y$ is the set of $g\in G$ such that $gG_eg^{-1}\sim G_e$.  If $\calc$ is the equivalence class
associated to $Y$, then $G_Y$ is the stabilizer of $\calc$ for the action of $G$  by conjugation on $\cale/\sim$. 

In Examples 
\ref{sec_deuxbouts} and \ref{sec_commens}, the group $G_Y$ is the commensurator of $G_e$, for any $e\subset Y$. 
In Examples
\ref{sec_hyp_rel} and \ref{sec_hyp_abs},  it is the maximal elementary subgroup containing $G_e$. 
In Example \ref{sec_ct}, it is the normalizer of the maximal abelian subgroup $A$ containing $G_e$ (it equals $A$ if $G$
is CSA). In Example \ref{sec_fini}, it is the normalizer of $G_e$.

We first note:
\begin{lem}\label{fact_fini}
  If $T$ is minimal, then $G_Y$ acts on $Y$ with finitely many orbits of edges.
\end{lem}

\begin{proof}
By minimality, there are finitely many $G$-orbits of edges in $T$.
If two edges of $Y$ are in the same orbit under some $g\in G$, then $g$ preserves $Y$, so they are in the same orbit 
under $G_Y$. 
\end{proof}

\begin{prop}\label{prop_acyl}
  Given $T$, the following statements are equivalent:
\begin{enumerate}
\item $T_c$ belongs to the same deformation space as $T$; 
\item  every stabilizer $G_Y$  is
elliptic in
$T$;
\item every cylinder  $Y\inc T$ is bounded;
\item no cylinder contains the axis of a hyperbolic element of $G$.
\end{enumerate}
\end{prop}

\begin{proof}
 We have seen $1\Leftrightarrow 2$.
If $Y$ is bounded, then $G_Y$ fixes a point of $T$ (the ``center'' of $Y$). If $G_Y$ fixes a point, then $Y$ is 
bounded by Lemma
\ref{fact_fini}. This shows
$2\Leftrightarrow 3$.

The implication $3\Rightarrow 4$ is clear. For the converse, assume that $Y$ is an unbounded cylinder. 
We know that $G_Y$ does
not fix a point. If all its elements are elliptic, then $G_Y$ fixes an end of $Y$. Any ray going out to that end maps
injectively to
$Y/G_Y$, contradicting Lemma
\ref{fact_fini}. Thus
$G_Y$   contains a hyperbolic element
$g$.  Being
$g$-invariant, $Y$ contains the axis of $g$.
\end{proof}

\begin{prop}\label{prop_petits}
Assume that any two groups of $\cale$ whose intersection is infinite are equivalent. Let $H$ be a
subgroup of  $  G$ which   is not virtually cyclic and is not an infinite, locally finite, torsion group. 

If $H$ is small, or commensurates an infinite small subgroup $H_0$, then $H$ fixes a point  in $T_c$.
 \end{prop}
 
 The hypothesis on $\cale$ is satisfied in Examples \ref{sec_deuxbouts} to \ref{sec_fini}, with the exception of
\ref{sec_commens}. Conversely, we will see in  Section \ref{sec_exa} that, in many examples,  groups which are elliptic in $T_c$ but not in $T$ are small.

Besides small groups, the proposition applies to groups $H$ which act on locally finite trees with small infinite stabilizers, for instance generalized Baumslag-Solitar groups.  It also applies to groups with a  small infinite normal subgroup.

\begin{proof} The result is clear if $H$ fixes a point of $T$.  
If  not, we show that $H$ preserves a subtree $Y_0$ contained in a cylinder.
This cylinder will  be $H$-invariant.

 First suppose that $H$ is small.  If $H$ preserves a line $\ell$,
  then $G_e\cap H$ is the same for all edges in that line. If $G_e\cap
  H$ is infinite, then $\ell$ is contained in a cylinder and one can take $Y_0=\ell$.
If $G_e\cap H$ is finite, then $H$ is virtually cyclic, which is ruled out.  By smallness, the only
  remaining possibility is that $H$ fixes a unique end of $T$.

  If $G_e\cap H$ is infinite for some edge $e$, we let $\rho $ be the
  ray joining $e$ to the fixed end. 
The assumption on $\cale$ implies that $\rho$ is contained in a cylinder $Y_0$. 
This cylinder   is $H$-invariant since $h\rho\cap\rho$ is a ray for any $h\in H$. 

 If all groups $G_e\cap H$ are finite, there are two
  cases. If $H$ contains a hyperbolic element $h$, the action of $H$
  on its minimal subtree $T_H$ is an ascending HNN extension with
  finite edge groups.  It follows that $T_H$ is a line  and $H$ is virtually cyclic.
  If every element of $H$ is elliptic,
  consider any finitely generated subgroup $H_0\subset H$. It fixes
  both an end and a point, so it fixes an edge. We conclude that $H_0$
  is finite, so $H$ is locally finite. This is ruled out.

  Now suppose that $H$ commensurates a small subgroup $H_0$. If $H_0$
  preserves a  unique line or fixes a unique end, the same is true for $H$ and
  we argue as before. If $H_0$ fixes a point $x\in T$, let $Y_0$ be
  the convex hull of the orbit $H.x$. Any segment $I\inc Y_0$ is
  contained in a segment $[hx,h'x]$ with $h,h'\in H$, and its
  stabilizer contains $hH_0h\m\cap h'H_0h'{}\m$ which is commensurable
  to $H_0$ hence infinite.  The assumption on $\cale$ implies   that $Y_0$ is contained in a cylinder.
\end{proof}

\begin{rem} 
Assume that there exists $C$ such that any two groups of $\cale$ whose intersection has order $>C$ are equivalent. 
The same proof shows that  locally finite subgroups are elliptic in $T_c$. 
\end{rem}

\subsection{The   collapsed tree   of cylinders $T\cc$}

Recall that   the collapsed tree of cylinders $T\cc$ is  the tree obtained from $T_c$ by collapsing all edges whose stabilizer is not in
$\cale$    (Definition \ref{dfn_T0}).

\begin{prop}\label{prop_bor} 
Cylinders of $T\cc$ have diameter at most 2.
\end{prop}

\begin{proof} 
 This follows from Remark \ref{rem_classe}:  if in $\cale$, the stabilizer of an edge $ (x,Y)$ of $T_c$  belongs to
  the equivalence class $\calc$ associated to $Y$.
  \end{proof} 

We say that    $\cale$ is \emph{sandwich-closed} if $A\subset H\subset B$ with $A,B\in\cale$ implies $H\in \cale$.
All families considered in Examples \ref{sec_deuxbouts} through \ref{sec_fini} have this property.

Sandwich-closedness has the following consequence. If $\varepsilon$ is an edge of $T_c$ such that  $G_\varepsilon$ is contained in a group of $\cale$, then $G_\varepsilon\in\cale$. This follows from Remark \ref{rem_gros}, which asserts that  $G_\varepsilon$ contains a group of $\cale$.

\begin{lem} \label{lem_fonct}
Assume that $\cale$ is  sandwich-closed. Given an equivariant map $f:T\to T'$, the cellular map $f_c:T_c\to T'_c$ of Proposition \ref{prop_cy} induces a cellular map $f\cc:T\cc\to T_c'^*$. 
\end{lem}

\begin{proof}  If $\varepsilon$ is  an edge of $T_c$ which is collapsed in $T\cc$,   
its image by $f_c$ is a point or an edge $\varepsilon'$ with $G_\varepsilon\inc G_{\varepsilon'}$. 
The group $G_\varepsilon$ is not in $\cale$, but by Remark \ref{rem_gros} it contains an element of $\cale$. 
Sandwich-closedness implies $G_{\varepsilon'}\notin\cale$, so $\varepsilon'$ is collapsed in $T_c'^*$. 
This shows that the natural map $T_c\to T_c'^*$ factors through $T\cc$.
\end{proof}

A subtree $X\subset T$ of diameter exactly 2 has a \emph{center}   $v\in V(T)$, 
and all its edges contain $v$. We say that $X$ is \emph{complete} if it contains all edges around $v$, \emph{incomplete} otherwise.

\begin{prop} Assume that $\cale$ is sandwich-closed. Let $T$ be  a minimal  $\cale$-tree.  
 \begin{enumerate*}
\item Every cylinder of $T\cc$ has diameter exactly $2$.
No stabilizer of an incomplete cylinder of $T\cc$   lies in $\cale$.
\item Conversely, assume that all cylinders of $T$ have diameter exactly $2$, 
and $G_Y\notin \cale$ for all incomplete cylinders $Y\subset T$.
Then $T\cc=T$.
  \end{enumerate*}
\end{prop}

\begin{proof} 
  It follows from Remark \ref{rem_classe} that any cylinder $Z$ of $T\cc$ has diameter at most $2$,
  and is obtained   from the ball of radius one around some $Y\in V_1(T_c)$   by collapsing all   edges
  with stabilizer outside $\cale$. It is incomplete if and only if at least one edge is collapsed. Note that the cylinders $Y\inc T$ and $Z\inc T\cc$ have the same stabilizer $G_Y=G_Z$. 
  
  We show that $Z$ has diameter exactly 2. Otherwise
   $Z$ consists of a single edge. 
The corresponding edge  $\eps=(v,Y)$ of $T_c$ is  
  the unique edge incident on $Y$ with $G_\eps\in\cale$, so $G_\eps=G_Y$ (otherwise, one would obtain other edges by applying elements of $G_Y\setminus G_\varepsilon$). By minimality of $T_c$, there exist other edges
  $\eps'$ incident on $Y$. They satisfy  $G_{\eps'}\notin\cale$, and   $G_{\eps'}\subset G_Y\in\cale$ contradicts sandwich-closedness.

If $Z$ is an incomplete cylinder of $T\cc$, at least one edge $\eps$ of $T_c$ incident on $Y$ is collapsed in $T\cc$,
 so $G_{\eps}\notin\cale$. As above, sandwich-closedness  implies $G_Y\notin \cale$. This proves 1.

 To prove 2, we shall define an isomorphism $g:T\cc\to T$. We denote by $v_Y$ the center of a cylinder $Y$ of $T$. 
  Let $f:T_c\ra T$ be the map sending $Y\in V_1(T_c)$ to $v_Y\in T$, sending $v\in V_0(T_c)$ to $v\in T$, and mapping the 
  edge $\varepsilon=(v,Y)$ to $[v,v_Y]$. Note that $ [v,v_Y]$ is an edge if $v\neq v_Y$, and is reduced to a point otherwise.
  
  We first prove that an edge $\varepsilon$ of $T_c$ is collapsed by $f$ (\ie $v=v_Y$)  if and only if $G_\varepsilon\notin\cale$. If $G_\varepsilon\notin \cale$, sandwich-closedness implies that $\varepsilon$ is collapsed, since   $T$ is an $\cale$-tree.
Conversely, if $\eps=(v,Y)$ is collapsed, then $v_Y=v\in V_0(T)$, so $v_Y$ lies in several cylinders of $T$. This implies that   $Y$ is incomplete, so $G_Y\notin\cale$. Since $G_Y$ is contained in $G_{v_Y}=G_v$, one has 
  $G_\eps=G_Y\notin\cale$.

  It follows that $f$ factors through a map $g:T\cc\ra T$ which maps edge to edge (without collapse). By minimality, $g$ is onto. 
There remains to prove that $g$ does not fold. If two edges of $T\cc$ have the same  image, they belong to the same cylinder. But $g$ is injective on each cylinder since $(v,Y)$ is mapped to $[v,v_Y]$.
\end{proof}

\begin{cor}Let $\cale$ be sandwich-closed. For any minimal $\cale$-tree $T$, one has $(T\cc)\cc=T\cc$. \qed
\end{cor}

\begin{prop}\label{prop_zer}
 Assume that $G_Y$ fixes a point of $T$ whenever there is an edge $\varepsilon=(x,Y)$ of
$T_c$ whose stabilizer is not in
$\cale$. Then $T_c$ and $T\cc$ belong to the same deformation space.
  Moreover, given $Y$, at most one edge   $\eps=(x,Y)$ of $T_c$ 
is collapsed in $T\cc$;  it satisfies $G_\eps=G_Y$. \end{prop}

\begin{proof} 
 Let $\varepsilon =(x,Y)$ be an edge of $T_c$ such that $G_\eps\notin\cale$. 
 It suffices to prove that $G_\eps=G_Y$ and that $G_{\eps'}\in \cale$
for every   edge $\eps'=(x',Y)$ with $x'\ne x$. 

By assumption, $G_Y$ fixes a point in $T$, hence a point $v\in Y$. 
If $x\ne v$, let $e$ be the initial edge of the segment $[x,v]$. 
Then $G_\varepsilon= G_x\cap G_Y\subset G_x\cap G_v\subset G_e$. On
the other hand, $G_e$ fixes $x$ and leaves $Y$ invariant, so $G_e\subset G_\varepsilon $. We conclude $G_\varepsilon
=G_e\in\cale$, a contradiction.
Thus $x=v$, and $\eps=(x,Y)$ is the only edge  incident to $Y$ with stabilizer not in $\cale$.
Moreover, since $G_Y$ fixes $x$, we have $G_{\varepsilon}=G_Y\cap G_x=G_Y$.
\end{proof}

\begin{cor}\label{cor_zerr}
If $T_c$ is in the same deformation space as $T$, then so is $T\cc$, and  therefore  $(T\cc)_c=T_c$ by Corollary \ref{cor_deformation}. \qed
 \end{cor}
 
 \begin{rem}   Suppose that $\cale$ is sandwich-closed and that, for any $A\in\cale$, any subgroup containing $A$ with index 2 lies in $\cale$. 
Then the hypothesis of Proposition \ref{prop_zer} is always satisfied when $G_Y$ is small. 
To see this, we suppose that $G_Y$ is not elliptic in $T$ and we show $G_\varepsilon\in\cale$. 
If $G_Y$ fixes an end of $T$, its subgroup $G_\varepsilon$, being elliptic, fixes an edge. 
  Since $G_\cale$ contains a group in $\cale$, sandwich-closedness implies
$G_\varepsilon\in\cale$. 
If $G_Y$ acts dihedrally on a line, some subgroup of $G_\varepsilon$ of index at most 2 fixes an edge, so  $G_\eps\in \cale$. 
\end{rem}

 Recall that cylinders of $T\cc$ have diameter at most $2$. We show that $T\cc$ is maximal for this property.

\begin{prop}\label{prop_maxi}
Assume that  $\cale$ is sandwich-closed. 
If $T'$ is any   $\cale$-tree dominated by $T$, and    cylinders of $T'$ are bounded, then $T'$ is  dominated by
  $T\cc$.  
 \end{prop}

\begin{proof}   
By Proposition \ref{prop_acyl} and  Corollary \ref{cor_zerr},  the tree  $T_c'^*$ belongs to the same
deformation space as $T'$.   Lemma \ref{lem_fonct} shows that $T\cc$ dominates $T_c'^*$, hence $T'$. 
\end{proof}

\subsection{Acylindricity}

 We now consider 
 acylindricity in the sense of Sela. 
Recall \cite{Sela_acylindrical} that  a tree is \emph{$k$-acylindrical} if the stabilizer of any segment of length $>k$ is trivial. It is
acylindrical if it is $k$-acylindrical for some $k$.  To handle groups with torsion, we  say that $T$ is 
 \emph{almost $k$-acylindrical} if  the stabilizer of any segment of length $>k$ is finite.

\begin{prop}\label{prop_acy}
 Assume that any two groups of $\cale$ whose intersection is  infinite
are equivalent.  Let $T$ be any $\cale$-tree.
\begin{enumerate}
\item The tree   $T\cc$  is almost 2-acylindrical. 
\item  If cylinders of $T$ are bounded (resp.\ of diameter $\leq k$), then $T$ is almost acylindrical (resp.\ almost $k$-acylindrical).
\end{enumerate}
 \end{prop}
 
Recall that the hypothesis on $\cale$ is satisfied in Examples \ref{sec_deuxbouts} to \ref{sec_fini}, with the exception
of \ref{sec_commens}.   The next section will provide examples where the converse to Assertion 2 holds.

 \begin{proof}

The first statement follows from the second one and Proposition  \ref{prop_bor}.

Since there are only finitely many orbits of cylinders, consider $k$ such that
cylinders have diameter at most $k$.
Any segment $I$ of length $k+1$ contains edges in distinct cylinders.
By the assumption on $\cale$,  the
 stabilizer of $I$ is finite. 
\end{proof}

 We also note:
\begin{lem} \label{lem_small}
Let $H$ be small, not virtually cyclic, not locally finite. Then $H$ is elliptic in any   almost acylindrical tree  $T$. 
\end{lem}

\begin{proof} The hypotheses on $H$ are the same as in Proposition \ref{prop_petits}, and the proof is similar. If $H$ is not elliptic in a tree $T$, it acts on a line with infinite edge stabilizers, or it fixes a unique end and some edge stabilizer is infinite. Both are  impossible if $T$ is almost acylindrical.
\end{proof}

\section{Examples} \label{sec_exa}

We now study specific examples. In most cases, we show  that $T\cc$ is  equal to $T_c$ (or at least in the same deformation space),
and we describe how far the deformation space of $T_c$ is from that of $T$. 

Recall that $T$ always dominates $T_c$. They are in the same deformation space if and only if, for every cylinder $Y$, the group
$G_Y$ is elliptic in $T$. Note that, if
all groups in $\cale$ are infinite, any virtually cyclic $G_Y$ is elliptic in $T$ because  by Remark \ref{rem_gros} it contains some $G_e$ (with finite
index). 

We also show that $T\cc$, which is almost 2-acylindrical by Proposition \ref{prop_acy}, is maximal for this property: it dominates
any almost acylindrical tree which is dominated by $T$. This is because, in the examples, groups which are elliptic in $T\cc$ but
not in $T$ are small,  and Lemma \ref{lem_small} applies. Describing the deformation space of $T\cc$ may thus be interpreted as finding which subgroups must be made elliptic in order to make $T$ almost acylindrical.
One may ask in general whether    a maximal almost acylindrical tree dominated by a given $T$ always exists.

\subsection{Relatively hyperbolic groups}\label{sec_exrel}

  \begin{prop}\label{prop_acyl_hyp}
Let $G$ be   hyperbolic relative  to $H_1,\dots,H_n$. Let  
$\sim$ be  co-elementarity,  as in Example \ref{sec_hyp_rel}. Let $T$ be a   tree with infinite
elementary edge stabilizers, such that    each $H_i$  is elliptic in $T$.   Then: \begin{enumerate}
\item
Edge stabilizers of $T_c$ are infinite
elementary, so $T\cc=T_c$. 
\item   
$T_c$ belongs to the same deformation space as $T$. In particular, it has the same non-elementary vertex stabilizers as $T$.
 \item $T_c$ is almost $2$-acylindrical  (and dominates any almost acylindrical tree  which is dominated by $T$). 
\end{enumerate}
\end{prop}

\begin{proof} 
Let $\varepsilon =(x,Y)$ be an edge of $T_c$. Here $G_Y$ is the maximal elementary subgroup containing $G_e$, for any edge $e$   of $ Y$.  This shows that $G_\varepsilon$ is elementary. It is infinite because it contains an element of $\cale$ 
(Remark  \ref{rem_gros}). 

To prove 2, we must show that every $G_Y$ is elliptic in $T$. 
If parabolic, $G_Y$ is elliptic  by assumption. If virtually cyclic, it is elliptic  by a remark made above (it contains    some $G_e$ with finite index).
Assertion 2 follows (its second half  is a general fact \cite[Corollary 4.4]{GL1}).

Assertion 3  now follows from Proposition   \ref{prop_acy}  (the parenthesized statement is trivial in this case).
\end{proof}

If we do not assume that $H_i$ is elliptic in $T$, we get:

  \begin{prop}\label{prop_acyl_hypa}
Let $G$ be   hyperbolic  relative   to finitely generated  one-ended  subgroups $H_1,\dots,H_n$. Let  
$\sim$ be  co-elementarity,  as in Example \ref{sec_hyp_abs}. Let $T$ be a   tree with infinite
elementary edge stabilizers.   Then: \begin{enumerate}
\item
Edge stabilizers of $T_c$ are infinite
elementary,  so $T\cc=T_c$. 
\item   
 $T$ and $T_c$ have the same non-elementary vertex stabilizers.  A subgroup is elliptic in $T_c$ if and only if it is  elliptic in $T$, or  parabolic. 
In particular,  $T_c$ is in the same deformation space as $T$ if and only if every parabolic subgroup is elliptic in $T$.
 \item  $T_c$ is almost $2$-acylindrical.  If the $H_i$'s are small, then $T_c$ dominates any almost acylindrical tree $T'$ which is dominated by $T$. 
\end{enumerate}
\end{prop}

\begin{proof} It is still true that every $G_Y$ is a maximal elementary subgroup (so 1 holds), but a parabolic $G_Y$ may now fail to
  be elliptic in $T$.  As pointed out at the beginning of the section, every virtually cyclic $G_Y$ is elliptic in $T$.

  If $G_v$ is a non-elementary vertex stabilizer of $T$, then $v$ belongs to two cylinders (otherwise  $G_v$  would be contained in
  some $G_Y$), so $G_v$ is a vertex stabilizer of $T_c$. The converse is clear since a non-elementary vertex stabilizer of $T_c$
  fixes a vertex of $V_0(T_c)$, so is a vertex stabilizer of $T$.

 To prove 2, there remains to show that each $H_i$ is elliptic in $T_c$. 
   If it is not elliptic in $T$, there is an edge $e$ with $G_e\cap H_i$ infinite (recall that $H_i$ is one-ended). 
   In particular, $G_e\sim H_i$. The associated equivalence class $\calc$ is invariant under conjugation by elements of $H_i$, so  
$H_i$ preserves the cylinder containing $e$  hence is elliptic in $T_c$.

  Acylindricity again follows from Proposition \ref{prop_acy}. The second part of 3 holds provided every $H_i$ is elliptic in
  $T'$, in particular if $H_i$ is small 
  by Lemma \ref{lem_small}. 
\end{proof}

\subsection{Abelian splittings of  CSA groups}\label{sec_excsa}

\begin{prop}\label{prop_acyl_CSA}
  Let $G$ be a torsion-free CSA group. 
Let $\cale$ (non-trivial abelian groups) and $\sim$ (commutation) be as  in Example \ref{sec_ct}. Let $T$ be an $\cale$-tree. Then:
\begin{enumerate}
\item Edge stabilizers of $T_c$ are non-trivial and abelian, so $T\cc=T_c$.
\item
 $T$ and $T_c$ have the same non-abelian vertex stabilizers.  A subgroup is elliptic in $T_c$ if and only if it is  elliptic in $T$ or is a non-cyclic abelian group.
In particular,  $T_c$ is in the same deformation space as $T$ if and only if every non-cyclic abelian subgroup of $G$ is elliptic in $T$.
\item $T_c$ is $2$-acylindrical and   dominates any   acylindrical $\cale$-tree  $T'$ which is dominated by $T$.
\end{enumerate}
\end{prop}

\begin{proof}  
If $Y\in V_1(T_c)$ is a cylinder, its stabilizer $G_Y$ is the set of $g\in G$ such that $gG_eg^{-1}$ commutes with $G_e$, for $e$ any edge of $Y$. By the CSA property, $G_Y$ is the maximal abelian subgroup containing $G_e$.  Conversely, a non-cyclic abelian subgroup acts on $T$ with non-trivial edge stabilizers and therefore leaves some cylinder invariant. 
As in the previous proof, non-abelian vertex stabilizers are the same for $T$ and $T_c$.

Assertions 1 and 2 follow from these observations. The vertex stabilizers of $T_c$ are those of $T$,
the non-cyclic maximal abelian subgroups which are not elliptic in $T$, and possibly cyclic subgroups which are elliptic in $T$. 

 The tree $T_c$ is 2-acylindrical by Proposition \ref{prop_acy}. 
  A group $H$ which is elliptic in $T_c$ but not in $T$ is abelian and non-cyclic, hence elliptic in $T'$ by Lemma \ref{lem_small}. Assertion 3 follows. 
 \end{proof}

We also note the following result, which gives a converse to the second assertion of Proposition \ref{prop_acy}:

\begin{prop}\label{prop_cycs} Let $G$, $\cale$, and $T$ be as above.
The following are equivalent:
\begin{enumerate}
\item every non-cyclic abelian subgroup is elliptic;  
\item   cylinders of  $T$ are bounded (equivalently, $T_c$ is in the same deformation space as $T$);
\item 
$T$ is   acylindrical;
\item no non-trivial element of $G$   fixes a line.
\end{enumerate}
\end{prop}

\begin{proof} We have just seen $1\Leftrightarrow2$.  Proposition \ref{prop_acy} gives $2\Rightarrow 3$, and obviously $3\Rightarrow 4$. To prove $4\Rightarrow 2$,
suppose that $Y$ is an unbounded cylinder. By Proposition \ref{prop_acyl}, it contains the axis   of a hyperbolic
element $g$. Let
$e$ be an edge contained in that axis, and $A$ the maximal abelian subgroup containing $G_e$. Since $gG_eg^{-1}$
commutes with $G_e$, the CSA property implies $g\in A$. Thus $G_e$ fixes the axis, contradicting 4.  
\end{proof}

\subsection{Cyclic splittings of commutative transitive groups}

The relation $\sim$ now is commensurability,  as in Example \ref{sec_deuxbouts}, so $G_Y$ is the commensurator of $G_e$ if $e$ is an edge of a cylinder $Y$.
 
For $s\ne0$, we denote by $BS(1,s)$   the   solvable Baumslag-Solitar group $ \langle
a,t\mid  tat^{-1}=a^s\rangle$. It is commutative transitive if and only if $s\ne -1$. Note that $BS(1,1)={\bbZ}^2$.

  \begin{prop}\label{prop_acyl_cyclic}
Let $G$ be torsion-free and commutative transitive. Let $\cale$ be the class of infinite cyclic subgroups of $G$,
and let $\sim$ be commensurability as in Example \ref{sec_deuxbouts}.  Let $T$ be an $\cale$-tree. Then:
\begin{enumerate}
\item  $T\cc$ and $T_c$ are in the same deformation space.
\item Every non-cyclic vertex stabilizer of $T$ is a vertex stabilizer of $T_c$ and $T\cc$,
and every other non-cyclic vertex stabilizer of $T\cc$ is isomorphic to some  $BS(1,s)$.  Every subgroup isomorphic to $BS(1,s)$ is elliptic in $T_c$ and $T\cc$. In particular, $T_c$ and $T\cc$ belong to the same deformation space as $T$ if and only if every subgroup of $G$ isomorphic to a
$BS(1,s)$ is elliptic in $T$. 
\item      $T\cc$ is   2-acylindrical and
  dominates any   acylindrical $\cale$-tree which is dominated
    by $T$.
    \end{enumerate}
  \end{prop}

\begin{proof}
 We note  the following algebraic facts, whose proof is left to the reader. Let $Z\subset H$ be an infinite cyclic
subgroup of a commutative transitive torsion-free group, and let $A$ be the centralizer of $Z$. Then $Z\subset
A\subset Comm(Z)\subset N(A)$. If $A=Z$, then $Z$ is malnormal.

Consider  a cylinder $Y\subset T$, and a vertex $v\in Y$.  
All edge stabilizers $G_e$, for $e\subset Y$, are commensurable, hence have the same centralizer
$A$ by commutative transitivity. By the previous remark, one has $A\subset G_Y\subset N(A)$ since $G_Y$ is the
commensurator of
$G_e$.

\begin{lem}
Assume that  $G_v\cap G_Y$ is non-cyclic. Then $G_Y$ fixes $v$, and only $v$.
 Moreover: 
 \begin{enumerate}
\item  if $Y$ is the only cylinder containing $v$, then   $G_Y=G_v$  and no edge of $T_c$ incident to the vertex $Y\in V_1(T_c)$ gets collapsed in $T\cc$;
\item if $v$ belongs to two cylinders, the edge $\eps=(v,Y)$ of $T_c$ is collapsed in $T\cc$ (the vertex $Y$ ``disappears'' in $T\cc$).
\end{enumerate}
\end{lem}

\begin{proof}
We first show that  $G_v\cap A$ is non-cyclic.
Assume   that $G_v\cap A$ is cyclic, necessarily infinite since it contains $G_e$ for $e$ an edge of $Y$
adjacent to $v$. By the initial note above, $G_v\cap A$ is malnormal in $G_v$, so $G_v\cap G_Y=G_v\cap A$ is
cyclic, a contradiction.

Since $G_v\cap A$ is non-cyclic, $v$ is its unique fixed point. 
It is also   the unique fixed point of $A$ (which centralizes $G_v\cap A$), and of
$G_Y\subset N(A)$.  

The ``moreover'' is clear: the only collapsible edge of $T_c$ incident to $Y$ is   $(v,Y)$, which exists if and only if $v$ belongs to two cylinders.
\end{proof}

By Proposition \ref{prop_zer}, the  lemma implies that $T_c$ and $T\cc$ belong to the same
deformation space. Moreover, any vertex stabilizer $H$ of $T\cc$ which is not a vertex stabilizer of $T$
equals $ G_Y$ for some cylinder $Y$ such that $G_v\cap G_Y$ is cyclic for every vertex $v\in Y$.
The group  $G_Y$ acts on $Y$ with  all   edge and vertex stabilizers   infinite cyclic.
Since it is commutative transitive, it is easy to see that $G_Y$  must be isomorphic to $\bbZ$  or
 a $BS(1,s)$. Conversely, any $BS(1,s)$ is elliptic in $T_c$ by Proposition \ref{prop_petits}.

 To prove Assertion 2, there remains to show that any non-cyclic vertex stabilizer $G_v$ of $T$ is a vertex
stabilizer of $T_c$ and $T\cc$.  This is clear if $v$ belongs to two cylinders. 
If it belongs to a unique cylinder $Y$, the  lemma tells us that 
$G_v=G_Y$ is a vertex stabilizer of $T_c$ and of $T\cc$.

 Asssertion 3 now follows from Proposition \ref{prop_acy}  and  Lemma \ref{lem_small}.
\end{proof}

\subsection{Commensurability}

In our last examples   $\sim$ is again commensurability, but we do not make assumptions on $G$, so   our  results are less precise.

\begin{prop}\label{prop_cycc}
Let $\cale$ be the set of two-ended subgroups of $G$, and $\sim$
    be commensurability, as in Example \ref{sec_deuxbouts}.
 Given an $\cale$-tree $T$, the following are equivalent:
\begin{enumerate}
\item cylinders of  $T$ are bounded (equivalently, $T_c$ is in the same deformation space as $T$);
\item the commensurator of each edge stabilizer  is elliptic in $T$;
\item 
$T$ is almost acylindrical;
\item no element of infinite order fixes a ray.
\end{enumerate}
\end{prop}

\begin{proof} 
$1\Leftrightarrow 2$ is clear because $G_Y=Comm(G_e)$ if $e\subset Y$.  $1\Rightarrow 3 $ follows from Proposition \ref{prop_acy}, and $3\Rightarrow 4$ is clear.

To prove  $4\Rightarrow 1$, assume that some cylinder is unbounded. By Proposition
\ref{prop_acyl}, it contains the axis $A_g$ of a hyperbolic element $g$. Let $e_0$ be an edge of $A_g$, let $e_i=g^i(e_0)$,
and $H_i=G_{e_i}$. If $H_i\subset H_{i\pm 1}$ for some $i$, then $H_i$ fixes a ray, contradicting 4. 
If not, we can find $h_0\in H_0$ and $h_2\in H_2$ not fixing $e_1$, and $h=h_0h_2$ is hyperbolic. 
As $H_0$ and $H_2$ are commensurable, there is a finite index subgroup
$H\subset H_0\cap H_2$ which is normal in both $H_0$ and $H_2$.
 Since $h$ normalizes $H$, the fixed point set of $H$
contains the axis of $h$, and 4 does not hold. 
\end{proof}

  \begin{cor}\label{cor_acyl_cyclic}
      For any $\cale$-tree $T$, the  tree  $T\cc$ is almost 2-acylindrical and   dominates any almost acylindrical $\cale$-tree which is dominated
    by $T$.
  \end{cor}

\begin{proof} 
This follows from Proposition \ref{prop_maxi}.
\end{proof}

In Example \ref{sec_commens}, we get:

\begin{prop}\label{prop_acyl_commens}
Let $\sim$ be the commensurability relation, with   $\cale$ as in Example \ref{sec_commens}.
Let $T $ be an $\cale$-tree. 
\begin{enumerate}
\item $T_c$ belongs to the same deformation space as $T$ if and only if there exists $k$ such that  any segment $I$ of
length
$>k$  contains an edge $e$ with the index
$[G_e:G_I]$ infinite. 
\item Assume that edge stabilizers are finitely generated. Then $T_c$ belongs to the same deformation space as $T$ if
and only if no group commensurable to an edge stabilizer fixes a ray.
\end{enumerate}
\end{prop}  

In particular, $T_c$ belongs to the same deformation space as $T$  when every group in $\cale$ is infinite  and
$T$ is almost acylindrical. But without further hypotheses on $\cale$  we cannot claim that $T\cc$ is almost acylindrical.

\begin{proof} 
1 follows from the fact that a segment $I$ is contained in a cylinder if and only if $[G_e:G_I]$ is finite for every $e\subset I$.  
The proof of 2 is fairly similar to that of Proposition \ref{prop_cycc}, and left to the reader. 
Finite generation of edge stabilizers is used to construct the normal subgroup of finite index  $H$.
\end{proof}

\section{JSJ splittings }
\label{sec_jsj}
\subsection{Generalities}

We review basic facts about JSJ splittings  and JSJ deformation spaces. See \cite{GL3} for details.

 In order to define JSJ splittings, one needs a family of edge groups which is closed under
  taking subgroups. 
Since $\cale$ does not have this property, we introduce the following substitute.

\begin{dfn}\label{dfn_substable}
  The family   
$\cale$ is \emph{substable} if, whenever  $G$ splits over a group $A$
contained in a group $B\in\cale$, then $A\in \cale$.
\end{dfn}

\begin{rem}
  When we work relative to a family of subgroups (like in Example \ref{sec_hyp_rel}),
the splitting of $G$ in the definition  should be relative to this family.
\end{rem}

\begin{example} \label{ex_substable}
In Examples \ref{sec_deuxbouts}, \ref{sec_hyp_abs}, \ref{sec_ct} (with $G$ torsion-free),
$\cale $ is substable if and only if
$G$ is one-ended.   In Example \ref{sec_commens}, $\cale $ is
substable if and only if $G$ does not split over a group having infinite index in a group of $\cale$.
In Example \ref{sec_hyp_rel}, we restrict to
relative splittings, and $\cale$ is substable if and only if $G$ is
one-ended relative to the $H_i$'s (\ie there is no non-trivial tree with
finite edge stabilizers in which every $H_i$ is elliptic).
\end{example}

We fix $\cale$ and $\sim$, with 
  $\cale$ substable. All trees are assumed to be $\cale$-trees. Strictly speaking, we   consider $\ol \cale$-trees, where $\ol\cale$  consists of all groups contained in a group of $\cale$. But substability guarantees that every $\ol\cale$-tree is an $\cale$-tree.

A   subgroup $H\subset G$ is \emph{universally elliptic}   if it is elliptic in every
 tree. A
 tree is universally elliptic if all its edge stabilizers are.

A  tree is a \emph{JSJ tree over $\cale$} if it is universally elliptic, and maximal for
this property: it dominates every universally elliptic tree. 
When   $\ol\cale$ consists of all groups with   a given property (e.g.\ abelian,
slender, elementary), we use the words abelian JSJ, slender JSJ, elementary JSJ...

JSJ trees always exist when $G$ is finitely presented,   and sometimes when $G$ is only finitely generated (in particular in the situations studied below).
They belong to the same deformation space, called
the \emph{JSJ deformation space over $\cale$}. If $T_J$ is a JSJ tree, and $T'$ is any    tree, there is a tree
$\hat T$ which refines $T_J$ and dominates
$T'$.

All these definitions and facts extend to the relative case: given  a 
  collection of 
 subgroups, one only considers
 trees   in which these subgroups are  elliptic. 

\subsection{QH-vertices} 
A vertex stabilizer of a JSJ tree is \emph{flexible} if it is not 
universally elliptic, and  does not belong to $\ol\cale$.  A key fact of JSJ theory is that flexible vertex stabilizers often have a very special form.

\begin{dfn}
A vertex stabilizer $G_v$ is a QH-subgroup (and $v$ is a QH-vertex) if there is an     exact sequence $1\to F\to G_v\xrightarrow{\pi}  \Sigma \to 1$,
where  
$\Sigma =\pi _1(S)$ is the fundamental group of a hyperbolic 2-orbifold with boundary.  
  Moreover, each incident edge stabilizer is conjugate to a   subgroup of a boundary subgroup 
$B\subset G_v$,   defined as the preimage under $\pi$ of  
$ \pi _1(C)$, with $C$ a component of
$\partial S$. 
\end{dfn}

In Subsection \ref{sec_comp2} we will need
a description of flexible vertex stabilizers in the following cases  (see \cite{GL3} for proofs). 
Assume that  $G$ is one-ended.

$\bullet$ $G$ is torsion-free and  CSA.  A   flexible vertex
stabilizer $G_v$ of an abelian JSJ tree   is a QH-subgroup, with $S$ a surface and $F$ trivial: $G_v$ is isomorphic to
$\Sigma =\pi _1(S )$, where $S$ is a compact surface. 

$\bullet$  $G$ is  hyperbolic relative to slender subgroups. 
A   flexible vertex
stabilizer
$G_v$ of a slender JSJ tree is a QH-subgroup, with $F$ finite. 

  In both cases, every incident edge stabilizer is conjugate to a finite index subgroup of a boundary subgroup.
Boundary subgroups are two-ended, and maximal among small  subgroups of
$G_v$.  Every boundary subgroup  contains an incident edge stabilizer.

\subsection{Canonical JSJ splittings}

We now use   trees of cylinders to  make JSJ splittings canonical  (\ie we get trees which are invariant under automorphisms).  See 
\cite{GL3} for a proof that JSJ splittings exist under the stated hypotheses, and  a discussion of their flexible vertices. 

\begin{thm}
Let $G$ be   hyperbolic relative to $H_1,\dots ,H_n$. 
Assume that $G$ is one-ended relative to $H_1,\dots, H_n$. 
There is an elementary (resp.\ virtually cyclic) JSJ tree relative to $H_1,\dots,H_n$ which is invariant under 
the subgroup of $\Out(G)$ preserving the conjugacy classes of the $H_i$'s. 
\end{thm}

\begin{rem}
When $H_i$ is not slender, this allows non-slender splittings. 
Still, one can describe flexible subgroups of this JSJ  tree as QH-subgroups \cite{GL3}.
 \end{rem}

\begin{proof}
First consider the case where $\cale$ is the family  of infinite elementary subgroups,  as in Example \ref{sec_hyp_rel}.
 It is substable because $G$ is one-ended relative to the $H_i$'s (see Example \ref{ex_substable}).
Let $T$ be an elementary JSJ tree relative to the $H_i$'s, and $T_c$ its tree of cylinders for co-elementarity.

By Proposition \ref{prop_acyl_hyp}, the tree $T_c$ has elementary edge stabilizers 
and lies in the  JSJ deformation space.
It is universally elliptic as its edge stabilizers are either   parabolic, 
or are virtually cyclic and contain an edge stabilizer of $T$ with finite index (Remark  \ref{rem_gros}).
It is invariant under the subgroup of $\Out(G)$ preserving the conjugacy classes of the $H_i$'s
because the JSJ deformation space is. The theorem follows.

Now turn to the case where $\cale$ is the class of infinite cyclic subgroups, still substable because of one-endedness.   We start with a virtually cyclic JSJ tree (relative to the $H_i$'s), we let 
$T_c$ be its tree of cylinders  (for co-elementarity, not commensurability), and we consider   $T\cc$ obtained
  by collapsing edges of $T_c$ whose stabilizer is not virtually cyclic.
By Propositions \ref{prop_zer} and \ref{prop_acyl_hyp}, the trees 
$T\cc$, $T_c$ and $T$ lie in the same deformation space. Moreover, $T\cc$ is universally elliptic
because its edge stabilizers contain an edge stabilizer of $T$ with finite index (Remark  \ref{rem_gros}).
It follows that $T\cc$ is a canonical JSJ splitting.
\end{proof}

 A similar argument, using  Proposition \ref{prop_acyl_CSA}, shows:

\begin{thm}
  Let $G$ be a one-ended torsion-free CSA group. There exists an abelian (resp.\ cyclic) JSJ tree of $G$ relative to all non-cyclic abelian subgroups, which is  $\Out(G)$-invariant. \qed
  \end{thm}

\section{Compatibility}\label{sec_compat}

Recall that two  trees $T$ and $T'$ are \emph{compatible} if they have a common refinement. 
The goal of this section is to show that $T_c$ is compatible with many splittings.
In particular, we show that   trees of cylinders of   JSJ deformation spaces     often are
\emph{universally compatible}, that is   compatible with every $\cale$-tree.

This is proved under two different types of hypotheses: 
in Subsection \ref{sec_comp1}  we assume that $\sim$ preserves universal ellipticity (this is true in particular when $\sim$ is commensurability),
and in Subsection \ref{sec_comp2} we work  in the setting of CSA groups and relative hyperbolic groups.  
Theorem \ref{thm_comp} follows from Corollary \ref{thm_compat}  and Theorem \ref{thm_JSJ_UC2} 

\subsection{A general compatibility statement}

We first prove a general compatibility statement, independent of JSJ theory. 
  \begin{prop}\label{prop_compat1} Let $T,T'$ be minimal
    $\cale$-trees.  If $T$ dominates $T'$, then $T_c$ is compatible
    with $T'$  and $T'_c$.
  \end{prop}

\begin{proof}  We have to construct a common refinement $\hat T$ of $T_c$ and $T'$ 
(compatibility of $T_c$ with $T'_c$ will follow since $T$ dominates $T'_c$). 
  Choose a map $f:T\to T'$ as in Subsection \ref{sec_fonct}.  For each $p\in V(T_c)$, denote by $Y_p$ the following
subset of
$T$:  the point $p$   if
$p\in V_0(T_c)$,   the cylinder  defining $p$  if
$p\in V_1(T_c)$. Consider $Z_p=f(Y_p)\subset T'$.
  By Lemma \ref{lem_image_cyl}, it is  either a point or a cylinder of $T'$. 
Note that a given edge of $T'$ is contained in exactly one  $Z_p$.

We obtain $\hat T$ from $T_c$ by ``blowing up'' each vertex $p$ to the subtree $Z_p$. Formally, we define $\hat T$ as the
tree obtained from $T_1=\Dunion_{p\in V(T_c)}Z_p$ as follows: for each edge 
$pq $ of $T_c$,  with $Y_p=\{x\}$ and $Y_q$ a cylinder containing $x$,  
add an edge to $  T_1$, the endpoints being attached to the two copies of $f(x)$ in $Z_p$ and $Z_q$ respectively.

The tree $T_c$ can be recovered from $\Hat T$ by collapsing each $Z_p$ to a point. We  show that $\hat T$ is also a
refinement of $T'$. Let $g:\Hat T\ra T'$ be the map defined as being induced by the identity on $T_1$, and being  
constant   on each added edge. It preserves alignment because a given edge of $T'$ is contained in  exactly one $Z_p$,
and $g $ is injective on each $Z_p$. One therefore recovers $T'$ from $\hat T$ by collapsing the added edges. 
\end{proof}

\subsection{Universal compatibility when $\sim$ preserves universal ellipticity}\label{sec_comp1}

As before, we fix $\cale$ and $\sim$, with 
  $\cale$ substable. All trees are assumed to be $\cale$-trees.
In this subsection, we   assume  
that $\sim$ preserves universal ellipticity in the following
sense. 

\begin{dfn} 
The relation
$\sim$ \emph{preserves universal ellipticity} if, given $A, B\in \cale$ with $A\sim B$,
the group $A$ is  universally elliptic   if and only if $B$   is.
\end{dfn}

For instance, commensurability always preserves universal ellipticity.
In the case of  a relatively hyperbolic group $G$, co-elementarity preserves  universal
ellipticity  if one restricts to
trees in which each $H_i$ is elliptic (Example \ref{sec_hyp_rel}). 
Similarly, in Example \ref{sec_ct}, one has to restrict to trees in which non-cyclic
abelian subgroups are elliptic  (the next subsection will provide non-relative results).

\begin{prop}\label{prop_JSJ_UC_gene}  
Assume that $\cale$ is substable and $\sim$ preserves universal ellipticity. 
If $T_J$ is a JSJ tree over $\cale$, 
its tree of cylinders is compatible with any $\cale$-tree.
\end{prop}

Using Example \ref{ex_substable}, we immediately deduce:  

  \begin{cor}\label{thm_compat} Let $G$ be finitely presented.
    \begin{enumerate}
    
   \item   Let  $\cale$ be the class of two-ended subgroups as in Example \ref{sec_deuxbouts}. If $G$ is one-ended, 
   the  tree of cylinders of the JSJ
      deformation space over $\cale$ is compatible with any
    tree with two-ended edge stabilizers.

    \item     More generally, let $\cale$ and $\sim$ be as in Example \ref{sec_commens}.
If  $G$ does not split over a subgroup having infinite index in a group of $\cale$,
 the tree of cylinders of the JSJ
      deformation space over $\cale$ is compatible with any
      $\cale$-tree.

    \item Let $G$ be hyperbolic relative 
      to finitely generated subgroups $H_i$ as in
      Example \ref{sec_hyp_rel}, and assume that $G$ is one-ended
      relative to the $H_i$'s.
      Then the tree of cylinders of the elementary JSJ deformation
      space relative to the $H_i$'s is compatible with any $
      \cale$-tree in which each $H_i$ is elliptic. \qed

    \end{enumerate}
  \end{cor}

\begin{rem} Finite presentability of $G$ is required only to know that  the JSJ deformation space exists. \end{rem}

\begin{proof}[Proof of Proposition \ref{prop_JSJ_UC_gene}] 
Let   $T$ 
be an $\cale$-tree, and let $\hat T$ be a
refinement of
$T_J$ which dominates $T$. 
Let $X$ be the tree obtained from $\hat T$ by collapsing all the edges whose stabilizer is not $\sim$-equivalent to an
edge stabilizer of $T_J$. The collapse map from $\hat T$ to $T_J$ factors through the collapse map $p:\hat T\to  X$. In
particular,
$X $ dominates $T_J$.

Since $T_J$ is universally elliptic, and $\sim$ preserves universal ellipticity, $X$ is universally elliptic.
By maximality of the JSJ deformation space, $X$ lies in the JSJ deformation space.
 In particular, $X$ and $T_J$ have the same tree of cylinders $X_c$.
We have to show that $X_c$ is compatible with $T$.

Let $\hat T _c$ be the tree of cylinders of $\hat T$, which is
compatible with $T$ by Proposition \ref{prop_compat1} since $\hat T$ dominates $T$. 
Because of the way $X$ was defined, the restriction of   $p:\Hat T\ra X$ to any cylinder is either constant or
injective.  By  Remark \ref{lem_fc_collapse}, $p_{c}:\hat T _c\to  X _c$ is a collapse map, so $X_c$ is compatible with
$T$.
\end{proof}

\subsection{Universal compatibility when $\sim$ does not preserve universal ellipticity} \label{sec_comp2}

\begin{thm}\label{thm_JSJ_UC2} Let $G$ be one-ended.  
\begin{enumerate}
\item  
Suppose $G$ is hyperbolic relative to \emph{slender} subgroups $H_1,\dots, H_n$.
The tree of cylinders of the slender JSJ deformation space     is compatible with every
tree whose edge stabilizers are slender.
\item  
Suppose $G$ is torsion free and CSA. The tree of cylinders of the abelian JSJ
deformation space    is compatible with every tree whose edge stabilizers are abelian.
\end{enumerate}
\end{thm}

The tree of cylinders is defined with $\cale$ as in Examples  \ref{sec_hyp_abs} and
\ref{sec_ct}: it   
  consists of all infinite slender (resp.\ abelian) subgroups, and $\sim$ is co-elementarity (=co-slenderness) or
commutation (note that each $H_i$ is finitely-ended, so  $\cale$ is admissible by Lemma \ref{lem_arb}).
The family
$\cale$
 is substable because $G$ is one-ended, but
$\sim$   does not preserve universal ellipticity.
 
\begin{proof}  
Let $T_J$ be a JSJ tree, 
 and $T_c$ its tree of cylinders. If $x\in
V_0(T_c)$, we know that
$G_x\notin
\cale$ (see Subsection \ref{sec_alg}). On the other hand  $G_Y\in\cale$ if $Y\in V_1(T_c)$, and
 edge stabilizers of
$T_c$ belong to $\cale$  (see
Subsections \ref{sec_exrel} and \ref{sec_excsa}).

We now show that $T_c$ is universally elliptic. Let $\varepsilon =(x,Y)$ be an edge. Let $e\subset Y$ be an edge of $T_J$
adjacent to
$x$. We have $G_e\subset G_\varepsilon \subset G_x$.
If $G_x$
is universally elliptic, so is 
 $G_\eps$. 
Otherwise, $G_x$ is  flexible. It is associated to a 2-orbifold $S$ as described in Section
\ref{sec_jsj}, and   
$G_e$ has finite index in a boundary subgroup $B$.  Since $B$ is the unique maximal small subgroup of $G_x$ containing
$G_e$, it also contains $G_\varepsilon $. Thus $G_e$ has finite index in $G_\varepsilon $, and $G_\varepsilon $ is 
universally elliptic because $G_e$ is.

Given any $\cale$-tree $T$, we now construct a common refinement $\hat T$ of $T_c$ and   $T $ by blowing up $T_c$ as in
the proof of Proposition \ref{prop_compat1}. There are several steps.

\paragraph{Step 1.} We first  define a $G_p$-invariant subtree $Z_p\subset T$, for $p$ a
vertex of
$T_c$. 

If $p\in V_0(T_c)$, the group $G_p$ is not  in $\cale$.   Consider its action on $T$. It
fixes a unique point, or it is a  flexible vertex group of $T_J$ and has a minimal invariant subtree in $T$ (because it is
finitely generated). We define
$Z_p$ as that point or subtree.

If $p\in V_1(T_c)$, then $G_p$ belongs to $\cale$. The cylinder of $T_J$ defining $p$ corresponds to an equivalence class 
$\calc\in\cale/\sim$ (which contains
$G_p$) as in Subsection \ref{sec_alg}.  If this class corresponds to a cylinder of $T$ (\ie if there is an edge of $T$ with
stabilizer equivalent to $G_p$), we define
$Z_p$ as that cylinder. If not, we now show that $G_p$ fixes a   point of $T$; this point is necessarily unique (otherwise,
there would be a cylinder), and we take it as
$Z_p$. 

Recall that $G_p$ is in $\cale$, hence is abelian or slender. If it does not fix a point in $T$, it fixes an end or   
preserves a line. Furthermore, it contains a subgroup $G_e$, 
  for $e$ an edge of $T_J$.  This
subgroup is elliptic in $T$ because
$T_J$ is universally elliptic. Some subgroup of index at most 2 of $G_e$ fixes an edge of $T$, yielding a cylinder
associated to $\calc$, a contradiction. 

\paragraph{Step 2.} We now explain how to attach edges of $T_c$ to 
$T_1=\Dunion_{p\in V(T_c)}Z_p$.
Let
$\varepsilon =pq $ be an edge,  with $p\in V_0(T_c)$ and $q\in V_1(T_c)$. We show that $G_\varepsilon $ fixes a unique
point $x_\varepsilon $ in $Z_p$, and this point $x_\varepsilon $ belongs to  $Z_q$; we then attach the endpoints of
$\varepsilon
$ to the copies of $x_\varepsilon $ in $Z_p$ and $Z_q$. 

Note that $G_\varepsilon $ is elliptic in $T$ (because $T_c$ is universally elliptic), and preserves $Z_p$ and $Z_q$.
If $Z_p$ is not a point, then $G_p$ is flexible, so is  an extension $F\ra G_p\ra \Sigma$. As explained above,
$G_\varepsilon
$ is contained in a boundary subgroup $B_0\subset G_p$ with finite index.

We consider the action of $G_p$ on its minimal subtree $Z_p\subset T$. Every boundary subgroup $B\subset G_p$
contains some $G_e$ with finite index (with $e$ an edge of $T_J$), hence acts elliptically.   
Being normal and finite, the group
$F$ acts as the identity, so there is an induced action of $\Sigma $ on $Z_p$. 
For that action, boundary subgroups of $\Sigma $ are elliptic, and 
edge stabilizers   are
  finite or  two-ended because they are slender (resp.\ abelian). 
This implies that $B_0$, hence also $G_\varepsilon $, fixes a unique
point $x_\varepsilon $ of
$Z_p$ (see \cite[Theorem III.2.6]{MS_valuationsI}
for the case of surface groups; the extension to an orbifold group  is
straightforward, as it contains a surface group with finite index).

We now show $x_\varepsilon \in Z_q$. If not, $G_\varepsilon $ fixes the initial edge $e$ of the segment joining
$x_\varepsilon $ to its projection onto $Z_q$. The stabilizer of $e$ is equivalent to $G_q$, and $Z_q$ was defined as the
cylinder containing $e$, so it contains $x_\varepsilon $.  

\paragraph{Step 3.} We can now construct   $\Hat T$ by gluing edges of $T_c$ to $T_1$ as in the proof of Proposition
\ref{prop_compat1}. It refines
$T_c$, and there is a natural map $g:
 \Hat T\ra T$ which is constant on all the edges corresponding to the edges of $T_c$, and which is isometric in
restriction to each
$Z_p$. To show that it is a collapse map, it suffices to see that $Z_p$ and $Z_{p'}$ (viewed as subtrees of $T$)  cannot have
an edge $e$ in common if
$p, p'$ are distinct vertices of $T_c$. 

We assume they do, and we reach a contradiction. Let $e_p\subset \hat T$ be the copy of $e$ in $Z_p$.
Then $G_{e_p}\subset G_e$. One has $G_{e_p}\in\cale$ because $\cale$ is substable,  and   Axiom 2 implies
$G_e\sim G_{e_p}$. This also shows that $\hat T$ is an $\cale$-tree.

Note that     $p$ and $p'$ cannot both belong to $V_1(T_c)$, as $Z_p$ and $Z_{p'}$ then are points or distinct
cylinders of $T$. 
We may therefore assume   $p\in V_0(T_c)$.
Let $\eps $ be the initial edge
of the segment $[p,p'] \subset T_c$.
By connectedness of cylinders, the segment joining $e_p$ to $e_{p'}$ is contained in a cylinder of $\Hat T$.
Since this segment contains the edge of $\Hat T$ corresponding to $\eps$, we have $G_{\eps}\sim G_{e_p}$.
Thus   $\langle G_\eps ,G_{e_p}\rangle$ is a small subgroup of $G_p$, and therefore is contained with finite index in a
boundary subgroup $B$. By \cite{MS_valuationsI}, $G_{e_p}$ fixes a unique point of $Z_p$, contradicting the fact that
$G_{e_p}$ fixes $e_p$. 
\end{proof}

\def\cprime{$'$} \newcommand{\noopsort}[1]{}

\begin{flushleft}
Vincent Guirardel\\
Institut de Math\'ematiques de Toulouse\\
Universit\'e de Toulouse et CNRS (UMR 5219)\\
118 route de Narbonne\\
F-31062 Toulouse cedex 9\\
France.\\
\emph{e-mail:}\texttt{guirardel@math.ups-tlse.fr}\\[8mm]

Gilbert Levitt\\
Laboratoire de Math\'ematiques Nicolas Oresme\\
Universit\'e de Caen et CNRS (UMR 6139)\\
BP 5186\\
F-14032 Caen Cedex\\
France\\
\emph{e-mail:}\texttt{levitt@math.unicaen.fr}\\
\end{flushleft}

\end{document}